\documentclass{article}

\usepackage{graphicx}

\usepackage[T1]{fontenc}
\usepackage{microtype}
\usepackage{titlesec}

\usepackage[backend=biber, style=ieee]{biblatex}

\usepackage{amsmath,amssymb,amsthm}
\usepackage[left=2.5cm,top=2.3cm,right=2.5cm]{geometry}
\usepackage[UKenglish]{babel}
\usepackage{bbold}
\usepackage{enumerate}
\usepackage{hyperref}
\usepackage{physics}

\newtheorem{prop}{Proposition}

\usepackage{tikz}
\usetikzlibrary{cd}
\tikzcdset{arrow style=tikz
	, diagrams={>={Stealth[length=1.2mm, width=1mm]}}
}

\usepackage[normalem]{ulem}
\usepackage{cancel}

\usepackage{xcolor}
\definecolor{darkgreen}{rgb}{0, 0.6, 0}

\renewcommand{\real}{{\mathbb R}} %% real numbers

\newcommand{\del}{\partial}

\newcommand{\pa}[2]{\frac{\partial #1}{\partial #2}}
\newcommand{\diff}[2]{\frac{\text{d} #1}{\text{d} #2}}
\renewcommand{\var}[2]{\frac{\delta #1}{\delta #2}}
\newcommand{\dt}{\frac{\text{d}}{\text{d} t}}

\newcommand{\bracket}[2]{\left\{ #1, #2 \right\}}

\newcommand{\bx}{\mathbf{x}}

\newcommand{\bv}{\mathbf{v}}
\newcommand{\bu}{\mathbf{u}}

\newcommand{\bE}{\mathbf{E}}
\newcommand{\bB}{\mathbf{B}}
\newcommand{\bA}{\mathbf{A}}
\newcommand{\be}{\mathbf{e}}
\newcommand{\bb}{\mathbf{b}}

\newcommand{\bF}{\mathbf{F}}
\renewcommand{\bf}{\mathbf{f}}
\newcommand{\bL}{\mathbf{\Lambda}}

\newcommand{\bW}{\mathbf{W}}

\newcommand{\cF}{\mathcal{F}}
\newcommand{\cG}{\mathcal{G}}

\newcommand{\hcF}{\hat{\mathcal{F}}}
\newcommand{\hcG}{\hat{\mathcal{G}}}

\newcommand{\hbE}{\hat{\mathbf{E}}}
\newcommand{\hbB}{\hat{\mathbf{B}}}

\newcommand{\hbF}{\hat{\mathbf{F}}}

\newcommand{\bh}{{\boldsymbol{\eta}}}
\newcommand{\br}{{\boldsymbol{\rho}}}

\newcommand{\hnab}{\hat{\nabla}}

\newcommand{\fM}{\mathbb{M}}

\newcommand{\fC}{\mathbb{C}}

\newcommand{\fE}{\mathbb{E}}

\newcommand{\fL}{\mathbb{\Lambda}}

\newcommand{\fV}{\mathbb{V}}

\newcommand{\fA}{\mathbb{A}}
\newcommand{\fH}{\mathbb{H}}

\newcommand{\fW}{\mathbb{W}}

\newcommand{\tcurl}{\text{curl}}

\newcommand{\tdiv}{\text{div}}

\renewcommand{\d}{\;\text{d}}
\newcommand{\dx}{\;\text{d} \bx}
\renewcommand{\dv}{\;\text{d} \bv}

\newcommand{\vth}{v_\text{th}}

\newcommand{\Dt}{\Delta t \;}

\begin{document}

\title{The Linearized Vlasov--Maxwell System as a Hamiltonian System (\textit{Second Revision, Clean Version})}
\author{
	Dominik Bell$^{1,2}$, Martin Campos Pinto$^1$, Stefan Possanner$^1$, Eric Sonnendrücker$^{1,2}$ \\
	\\
	\normalsize{$^1$ Max-Planck-Institut für Plasmaphysik, Garching, Germany} \\
	\normalsize{$^2$ Technische Universität München, Zentrum Mathematik, Garching, Germany}\vspace{7mm}
	}

\maketitle

\begin{abstract}
	We present a Hamiltonian formulation for the linearized Vlasov--Maxwell system with a Maxwellian background distribution function. We discuss the geometric properties of the model at the continuous level, and how to discretize the model in the GEMPIC framework \cite{GEMPIC}. This method allows us to preserve the structure of the system at the semi-discrete level. To integrate the model in time, we employ a Poisson splitting and discuss how to integrate each subsystem separately.
	
	We test the model against the direct delta-f method, which is the non-geometric pendant of our model. The first test case is the weak Landau damping, where our model exhibits the same physical properties for short simulations, but enjoys better long-time stability and energy conservation due to its geometric construction. These advantages becomes even more pronounced for the simulation of Bernstein waves, our second test case, where the noise in the direct delta-f method washes out all features of the dispersion relation whereas our model is able to reproduce the full spectrum correctly. The model is implemented in the open-source \texttt{Python} library \texttt{STRUPHY} \cite{StruphyPaper, StruphyDoc}.
\end{abstract}

% !TeX spellcheck = en_GB

\section{Introduction}

Like many other physical models \cite{Morrison_2017}, the Vlasov--Maxwell system \cite{Vlasov_1968} possesses a rather rich structure, most prominently, there exists a non-canonical Hamiltonian system \cite{Marsden1982, Morrison_1982, Morrison_Ideal_Fluid} which consists of a Hamiltonian functional, that has the interpretation of physical energy, and a Poisson bracket, which generates the equations of motion. The Casimir functionals of this bracket correspond to physical conservation properties, most importantly the divergence-free magnetic field and Gauss' law.

Geometric methods aim to discretize the infinite-dimensional system so that the physical (Hamiltonian) structure is also present in the finite-dimensional approximation space, rather than simply discretizing the model equations of motion. This construction guarantees long-time stability due to the boundedness of the Hamiltonian \cite{Book_Geometric_Numerical_Integration}. Since their inception, geometric Particle-in-cell (PIC) methods have been discussed and extended manifold in the literature \cite{Xiao_2018, Markidis_2011, Squire_2012, Evstatiev_2013}. The structure-preserving discretization techniques used in this work are closely related to the GEMPIC \cite{GEMPIC, Perse_2021, Kormann_2021} framework. It consists of two components: the discretization of the distribution function $f$ using a PIC method and the Finite Element Exterior Calculus (FEEC) \cite{Arnold_Falk_Winther_2006, Arnold_Falk_Winther_2010} method using B-splines as basis functions \cite{Buffa_2011} for the discretization of the electromagnetic fields. Discretizing the Hamiltonian structure in this way guarantees structure-preservation at the semi-discrete level.

The PIC method is a popular choice for transport equations because markers are evolved along the characteristics of the equation and integrals are computed in the sense of a Monte Carlo expectation value. However, the distribution function $f$ appearing in the Vlasov--Maxwell model is defined on a six-dimensional phase space, which makes simulations very costly if one wants to achieve acceptable noise levels. There are (at least) two ways to improve upon this situation: either a model reduction such as averaging out the gyration of the particles around the magnetic field lines yielding 5-dimensional models \cite{Handbook_of_Plasma_Physics, Burby_2019} or considering the model near thermal equilibrium where a fluid description is a good approximation \cite{Kulsrud_1980}. However, this has the downfall of potentially not capturing all the physical effects. The other way is to improve the numerical algorithms. For particle methods, the goal is often to reduce the noise.

\begin{figure}[h]
	\centering
	\includegraphics[width=0.9\linewidth]{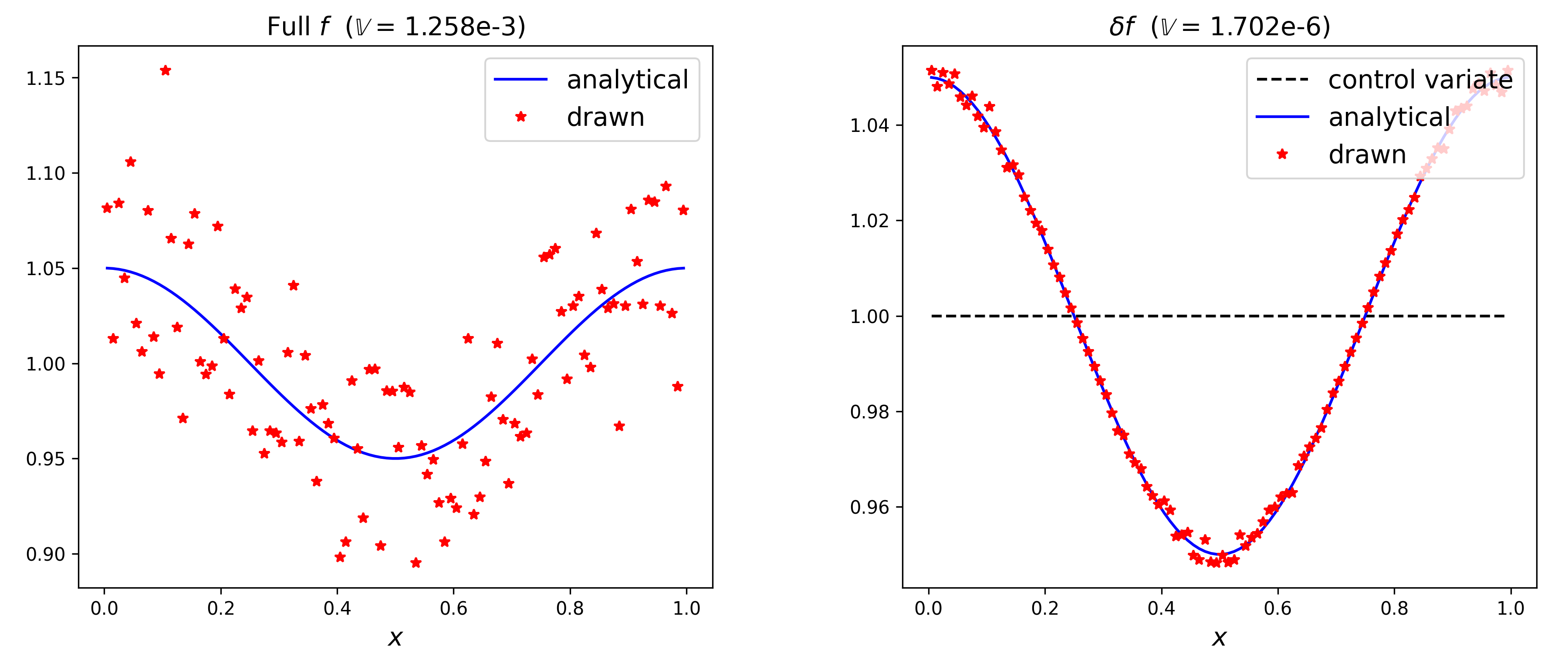}
	\caption{Sampled points from function \eqref{sampling-function} without (left) and with (right) using of the constant function $f_0(x)=1$ as a control variate. Drawing with a control variate lowers the variance $\fV$ with respect to the analytical function by a factor of $10^3$.}
	\label{fig:sampling}
\end{figure}

One such noise reduction method is the control variate method, which assumes that one can analytically evaluate expectation values against a function subtracted from the distribution function to compute integral values more accurately, thus increasing the overall accuracy. To illustrate the advantage of a control variate, let us look at an example: assume that we want to sample the function
\begin{equation}\label{sampling-function}
	f(x) = 1 + \delta \, \cos\left(2\pi \frac{x}{L}\right)
\end{equation}
on a domain $[0,L]$ with $\delta$ some small number, $\abs{\delta} \ll 1$. In this case, the control variate is the simple function $f_0(x) = 1$, and we only draw markers from $f_1(x) = \delta \, \cos\left(2\pi x / L\right)$. The difference between drawing from the full $f(x)$ and drawing only from $f_1(x)$ is represented in fig. \ref{fig:sampling} (both using the same pseudo-random sampling technique). For a number of markers $N_p = 10^5$ the control variate achieves a variance (with respect to the analytical solution) which is three orders of magnitude lower than the full-$f$ method. However, applying this procedure to the distribution function in the Vlasov--Maxwell model does \textit{a priori} break the structure of the model because the control variate terms are introduced into the equations.

The aim of this work is to bring the concepts of Hamiltonian structure-preservation and noise reduction using a control variate together. If it is known that the solution of the Vlasov-Maxwell system will stay close to a given equilibrium, one can consider the linearized model equations instead of solving the full system; this equilibrium is then used as the control variate. We present a Hamiltonian framework for the linearized Vlasov--Maxwell system using an equilibrium as a control-variate and discretize it in a structure-preserving manner. We compare our model to non-geometric discretizations of the non-linear Vlasov--Maxwell equations and of the linearized equations. For the weak Landau damping, we find that our model reproduces the correct physical predictions for short times, i.e., the damping rate and damping level are the same for both. However, the structure-preservation shows its advantages in terms of long-time stability and energy conservation. The weak Landau damping is simulated both in Cartesian and curvilinear coordinates. As a second test, we check for plasma waves appearing under a background magnetic field, where we see that our code reproduces correctly the dispersion relation including the electromagnetic waves as well as the Bernstein waves.

This work is organized as follows: In section \ref{sec:continuous}, we present the model, its Hamiltonian structure, and formulate both the equations of motion and the Hamiltonian system in curvilinear coordinates. Then we discuss the discretization of this geometric delta-f model in a structure-preserving way using the GEMPIC framework. The PIC method is presented in sec. \ref{sec:pic-discretization} where we also show how to handle the source term appearing in the linearized Vlasov equation. This source term causes the weights of the PIC-function to be time-dependent, adding another variable to the system compared to \cite{GEMPIC}. We derive the particle equations of motion, and subsequently discretize the field in sec. \ref{sec:feec-discretization}, yielding the fully semi-discrete system. In order to integrate the system in time with energy conservation, we employ a Poisson-splitting and show how to numerically integrate each subsystem them in time (sec. \ref{sec:poisson-splitting}). We also discuss implementation details as the code is publicly available in the open-source \texttt{Python}-library \texttt{STRUPHY} \cite{StruphyPaper, StruphyDoc}. In sec. \ref{sec:results}, we present numerical results for the weak Landau damping and the dispersion relation under a background magnetic field. Finally, we summarize our findings and give an outlook in sec. \eqref{sec:summary_outlook}.

% !TeX spellcheck = en_GB

\section{The Model}\label{sec:continuous}

\subsection{The Linearized Vlasov--Maxwell Model}

We study the linearized Vlasov--Maxwell model for one species on a physical domain $\Omega\subseteq\real^3$
\begin{subequations}\label{final-cart-maxwell-equations}
	\begin{align}
		\del_t f_1 + \bv \cdot \, & \nabla_\bx f_1 + \frac{1}{\epsilon} \left( \bE_0 + \bv\times\bB_0 \right) \cdot \nabla_\bv f_1 = \frac{1}{\vth^2 \, \epsilon} \bE_1 \cdot \bv f_0 \, , \label{final-cart-vlasov-equation} \\
		\del_t \bE_1 & = \nabla_\bx \times \bB_1 - \frac{\alpha^2}{\epsilon} \int_{\real^3} f_1 \, \bv \dv \, , \label{final-cart-ampere-equation} \\
		\del_t \bB_1 & = - \nabla_\bx \times \bE_1 \, , \label{final-cart-faraday-equation} \\
		\nabla_\bx \cdot \bE_1 & = \frac{\alpha^2}{\epsilon} \int_{\real^3} f_1 \dv \, , \label{final-cart-gauss-equation} \\
		\nabla_\bx \cdot \bB_1 & = 0 \, , \label{final-cart-magn-gauss-equation}
	\end{align}
\end{subequations}
such that the perturbations $f_1$, $\bE_1$, and $\bB_1$ are self-consistently coupled. The background functions $f_0, \bE_0, \bB_0$ are assumed to be an equilibrium solution of the Vlasov--Maxwell equations
\begin{align}\label{f0-is-maxwellian}
	f_0(\bx,\bv) & = \frac{n_0(\bx)}{\left( \sqrt{2\pi} \vth \right)^3} \exp\left( -\frac{|\bv|^2}{2\vth^2} \right) \, , & \bE_0(\bx) & = \vth^2 \, \epsilon \, \nabla_\bx \ln(n_0(\bx)) \, , & \bB_0 & = const. \, ,
\end{align}
such that we have in particular
\begin{equation}\label{equilibrium-condition}
	\bv \cdot \nabla_\bx f_0 + \frac{1}{\epsilon} \left( \bE_0 + \bv \times \bB_0 \right) \cdot \nabla_\bv f_0 = \bv \cdot \left( \frac{f_0}{n_0} \, \nabla_\bx n_0 - \frac{1}{\vth^2 \epsilon} \bE_0 \, f_0 \right) = 0 \, .
\end{equation}
The parameters capturing all the physical scales are
\begin{align}
	\frac{1}{\epsilon} & = 2\pi \frac{\Omega_c}{\hat{\omega}} \, , &\alpha^2 & = \frac{\Omega_p^2}{\Omega_c^2} \, , & \Omega_c & = \frac{q}{m} \hat{B} \, , & \Omega_p & = \sqrt{\frac{q^2 \hat{n}}{\epsilon_0 m}} \, ,
\end{align}
where hats denote a unit of each quantity, $q$ is the charge and $m$ the mass of the species, and the velocity scale is the speed of light which also determines the scales of $\bE$ in terms of $\hat{B}$, i.e.,
\begin{equation}
	\frac{\hat{E}}{\hat{B}} = c = \hat{v} \, .
\end{equation}

\subsection{Hamiltonian Structure}

The linearized Vlasov--Maxwell system conserves the following energy functional \cite{KruskalObermann1958, Shadwick_PhDThesis}
\begin{equation}\label{final-cart-Hamiltonian}
	H[f_1, \bE_1, \bB_1] = \frac{\alpha^2}{2} \vth^2 \int_{\Omega\times \real^3} \frac{f_1^2(t, \bx, \bv)}{f_0(\bx, \bv)} \dx \dv + \frac{1}{2} \int_\Omega \abs{\bE_1(\bx)}^2 \dx + \frac{1}{2} \int_\Omega \abs{\bB_1(\bx)}^2 \dx \, .
\end{equation}

\begin{prop}
	The following bracket is a Poisson bracket
	\begin{subequations}\label{final-cart-bracket}
		\begin{align}
			\bracket{\cF}{\cG} & = \int_\Omega \left[ \left( \nabla_\bx \times \var{\cF}{\bE_1} \right) \cdot \var{\cG}{\bB_1} - \var{\cF}{\bB_1} \cdot \left( \nabla_\bx \times \var{\cG}{\bE_1} \right) \right] \dx \\
			& \quad + \frac{1}{\vth^2 \, \epsilon} \int_{\Omega\times \real^3} \left[ \var{\cF}{f_1} \var{\cG}{\bE_1} \cdot \bv - \var{\cF}{\bE_1} \cdot \bv \var{\cG}{f_1} \right] f_0 \dx \dv \\
			& \quad + \frac{1}{2} \frac{1}{\alpha^2 \vth^2} \int_{\Omega\times \real^3} f_0 \left[ \bv\cdot \left( \left(\nabla_\bx \var{\cF}{f_1}\right) \var{\cG}{f_1} - \var{\cF}{f_1} \nabla_\bx \var{\cG}{f_1} \right) \right. \\
			& \hspace{35mm} + \frac{1}{\epsilon} \left. \left(\bE_0 + \bv\times \bB_0 \right) \cdot \left( \left(\nabla_\bv \var{\cF}{f_1}\right) \var{\cG}{f_1} - \var{\cF}{f_1} \nabla_\bv \var{\cG}{f_1} \right) \right] \dx \dv \, .
		\end{align}
	\end{subequations}
\end{prop}
Same as we suppressed arguments here, we will mostly leave out the integral domain since integrals are over the whole domain with the measure indicating which space is integrated over ($\dx$ for the spatial domain $\Omega$ and $\dx \dv$ for the whole phase-space $\Omega\times \real^3$).

\textit{Note:} The Hamiltonian \eqref{final-cart-Hamiltonian} and the Poisson bracket \eqref{final-cart-bracket} can be found by using the energy--Casimir method \cite{Holm1985, Morrison_Ideal_Fluid, Shadwick_PhDThesis}. In particular, by integration by parts, using \eqref{equilibrium-condition}, and symmetrizing the transport part, the Poisson bracket from the energy--Casimir method can be shown to be equivalent to \eqref{final-cart-bracket}. The form we mention is more suitable for discretization as it gives dynamical weights, which is also expected from having a non-zero right-hand side in a transport equation as discussed below.

\textit{Proof:} We prove that this bracket \eqref{final-cart-bracket} is indeed a Poisson bracket. It is manifestly antisymmetric and since it consists of functional derivatives in pairs it also verifies the Leibniz property. It remains to show that it satisfies the Jacobi identity
\begin{equation}
	\bracket{\bracket{U_i}{U_j}}{U_k} + \bracket{\bracket{U_k}{U_i}}{U_j} + \bracket{\bracket{U_j}{U_k}}{U_i} = 0
\end{equation}
for $U_i \in \left\{\bE_1, \bB_1, f_1\right\}$. Because the only non-vanishing bracket between two coordinates is
\begin{equation}
	\bracket{f_1}{E_{1,i}} = \frac{1}{\vth^2 \, \epsilon} \int v_i f_0 \d \bx \d \bv = - \bracket{E_{1,i}}{f_1} \, ,
\end{equation}
of which any further functional derivative vanishes
\begin{equation}
	\var{}{U_i} \bracket{f_1}{E_{1,j}} = 0 \, ,
\end{equation}
and so the bracket is indeed a Poisson bracket. \qed

\begin{prop}
	The above Hamiltonian \eqref{final-cart-Hamiltonian} and Poisson bracket \eqref{final-cart-bracket} form a non-canonical Hamiltonian system that generates the linearized Vlasov--Maxwell system consisting of \eqref{final-cart-vlasov-equation}, \eqref{final-cart-ampere-equation}, and \eqref{final-cart-faraday-equation} as its equations of motion while \eqref{final-cart-gauss-equation} and \eqref{final-cart-magn-gauss-equation} are Casimirs of the bracket.
\end{prop}

\textit{Proof:} It is a straightforward computation to see that the equations of motion are obtained in weak form by
\begin{equation}
	\dt \cF = \bracket{\cF}{H}
\end{equation}
for any functional $\cF = \cF[f_1, \bE_1, \bB_1]$.

\subsection{Curvilinear Coordinates}

Following \cite{Holderied2021}, we will formulate the model in curvilinear coordinates. We want to be able to solve the above system in a general class of geometries which can be described by a diffeomorphism $L$ from a logical domain $\hat{\Omega} = [0,1]^3$ with coordinates $\bh$ to the physical domain $\Omega$ with coordinates $\bx$:
\begin{equation}\label{mapping}
	L \, : \, \hat{\Omega} \to \Omega \,, \, \bh \mapsto L(\bh) = \bx \, .
\end{equation}
We denote its Jacobian matrix by $DL$, the metric matrix by $G = DL^T DL$, and its determinant $\sqrt{g} = \det(DL) = \sqrt{\det(G)}$.

In this coordinate transformation, we will treat the electromagnetic fields not as vector fields but as the coefficients of forms. This allows us to leverage the finite element exterior calculus (FEEC) framework, which guarantees the preservation of the geometric properties of our model upon discretization. We will discuss details below, and for now, just note that the fields live in an exact de Rham sequence
\begin{equation}\label{physical-de-Rham-sequence}
	\begin{tikzcd}
		H^1(\Omega) \arrow[r,"\;\nabla_\bx \;"] & H(\tcurl, \Omega) \arrow[r,"\;\nabla_\bx \times \;"] & H(\tdiv, \Omega) \arrow[r,"\;\nabla_\bx \cdot \;"] & L^2(\Omega) \, ,
	\end{tikzcd}
\end{equation}
which also holds after pull-back to the logical domain
\begin{equation}\label{logical-de-Rham-sequence}
	\begin{tikzcd}
		H^1(\hat{\Omega}) \arrow[r,"\;\hnab_\bh \;"] & H(\widehat{\tcurl}, \hat{\Omega}) \arrow[r,"\;\hnab_\bh \times \;"] & H(\widehat{\tdiv}, \hat{\Omega}) \arrow[r,"\;\hnab_\bh \cdot \;"] & L^2(\hat{\Omega}) \, ,
	\end{tikzcd}
\end{equation}
where we introduced the notation to denote the coefficients on the logical domain by the same symbols as the physical domain but with hats. Under a mapping such as \eqref{mapping}, coefficients of forms transform as follows:
\begin{subequations}
	\begin{align}
		L^\ast \phi & = \hat{\phi} = \phi && \forall \, \phi \leftrightarrow f^0 \in H^1(\Omega) \, , \\
		L^\ast \bu & = \hat{\bu} = DL^T \bu && \forall \, \bu \leftrightarrow u^1 \in H(\tcurl, \Omega) \, , \\
		L^\ast \bv & = \hat{\bv} = \sqrt{g} \, DL^{-1} \bu && \forall \, \bv \leftrightarrow v^2 \in H(\tdiv, \Omega) \, , \\
		L^\ast \phi & = \hat{\varphi} = \sqrt{g} \, \varphi && \forall \, \varphi \leftrightarrow p^3 \in L^2(\Omega) \, .
	\end{align}
\end{subequations}

We want to keep the magnetic Gauss law in strong form, so it follows
\begin{align}
	\bB_1 & \longleftrightarrow B^2 \in H(\tdiv) &&\text{and} & \bE_1 & \longleftrightarrow E^1 \in H(\tcurl) \, .
\end{align}
The background electromagnetic fields will be treated the same as their dynamical cousins, and the distribution function $f_1$ itself will be treated as a 0-form. We note that we cannot take the strong curl of the magnetic field, and so Ampère's equation \eqref{final-cart-ampere-equation} will have to be verified weakly.

On the logical domain, the dynamical Maxwell equations \eqref{final-cart-ampere-equation}, \eqref{final-cart-faraday-equation} then become
\begin{equation}\label{curvy-Faraday}
	\del_t \hbB_1 = - \hnab_\bh \times \hbE_1 \, ,
\end{equation}
and also $\forall \; \hbF \leftrightarrow L^\ast F^1 \in \hat{V}^1$ we have
\begin{equation}\label{curvy-Ampere}
		\dt \int_{\hat{\Omega}} \hbF \cdot G^{-1} \, \hbE_1 \, \sqrt{g} \d \bh = \int_{\hat{\Omega}} \left(\hnab_\bh \times \hbF \right) \cdot G\, \hbB_1 \frac{1}{\sqrt{g}} \d \bh - \frac{\alpha^2}{\epsilon} \int_{\hat{\Omega}\times \real^3} DL^{-T} \, \hbF \cdot \bv \, \hat{f}_1 \, \sqrt{g} \d \bh \dv \, .
\end{equation}
The modified Vlasov equation \eqref{final-cart-vlasov-equation} becomes
\begin{equation}\label{final-curvy-vlasov}
	\del_t \hat{f}_1 + \bv \cdot DL^{-T} \hnab_\bh \hat{f}_1 + \frac{1}{\epsilon} \left( DL^{-T} \hbE_0 + \bv \times \frac{1}{\sqrt{g}} \, DL \, \hbB_0 \right) \cdot \nabla_\bv \hat{f}_1 = \frac{1}{\vth^2 \, \epsilon} DL^{-T} \hbE_1 \cdot \bv \hat{f}_0 \, .
\end{equation}
The Hamiltonian \eqref{final-cart-Hamiltonian} becomes
\begin{equation}\label{final-curvy-Hamiltonian}
	\hat{H} = \frac{\alpha^2}{2} \vth^2 \int_{\hat{\Omega} \times \real^3} \frac{\hat{f}_1^2}{\hat{f}_0} \sqrt{g} \d\bh \dv + \frac{1}{2} \int_{\hat{\Omega}} \hbE_1 \cdot G^{-1} \hbE_1 \sqrt{g} \d\bh + \frac{1}{2} \int_{\hat{\Omega}} \hbB_1 \cdot G \hbB_1 \frac{1}{\sqrt{g}} \d\bh \, ,
\end{equation}
which is still preserved because in particular $\hnab_\bh \cdot \left( \sqrt{g} DL^{-1} \bv \right) = 0$.

Under the pull-back, functional derivatives transform dual to the variable to which the derivative is taken with respect to. The Poisson bracket \eqref{final-cart-bracket} is thus transformed to
\begin{subequations}\label{final-curvy-bracket}
	\begin{align}
			\bracket{\hcF}{\hcG} & = \int \left[ \left( \hnab_\bh \times \var{\hcF}{\hbE_1} \right) \cdot G \var{\hcG}{\hbB_1} - \left( \hnab_\bh \times \var{\hcG}{\hbE_1} \right) \cdot G \var{\hcF}{\hbE_1} \right] \frac{1}{\sqrt{g}} \d \bh \d \bv \\
			& \quad + \frac{1}{\vth^2 \, \epsilon} \int \left[ \var{\hcF}{\hat{f}_1} \left(DL^{-T} \var{\hcG}{\hbE_1}\right) \cdot \bv - \var{\hcG}{\hat{f}_1} \left(DL^{-T} \var{\hcF}{\hbE_1}\right) \cdot \bv \right] \hat{f}_0 \sqrt{g} \d \bh \d \bv\\
			& \quad + \frac{1}{2} \frac{1}{\alpha^2 \vth^2} \int \hat{f}_0 \left[ \bv\cdot DL^{-T} \left(\left(\hnab_\bh \var{\hcF}{\hat{f}_1}\right) \var{\hcG}{\hat{f}_1} - \var{\hcF}{\hat{f}_1} \hnab_\bh \var{\hcG}{\hat{f}_1} \right) \right. \\
			& \hspace{20mm} \left. + \frac{1}{\epsilon} \left(DL^{-T} \hbE_0 + \bv\times \sqrt{g} DL \, \hbB_0 \right) \cdot \left(\left(\nabla_\bv \var{\hcF}{\hat{f}_1}\right) \var{\hcG}{\hat{f}_1} - \var{\hcF}{\hat{f}_1} \nabla_\bv \var{\hcG}{\hat{f}_1} \right) \right] \sqrt{g} \d \bh \dv \, .
	\end{align}
\end{subequations}

It is straightforward to show that the equations \eqref{curvy-Faraday}, \eqref{curvy-Ampere}, and \eqref{final-curvy-vlasov} are the equations of motion of the Hamiltonian system formed by the Hamiltonian \eqref{final-curvy-Hamiltonian} and the bracket \eqref{final-curvy-bracket}.

% !TeX spellcheck = en_GB

\section{Particle Discretization}\label{sec:pic-discretization}

We will discretize the distribution function $f_1$ using particles. This immediately leaves us with two choices:
\begin{enumerate}
	\item Derive particle equations of motion from the discretized Vlasov equation \eqref{final-curvy-vlasov}
	\item Discretize the Hamiltonian structure (Hamiltonian functional \eqref{final-curvy-Hamiltonian} and Poisson bracket \eqref{final-curvy-bracket}) and derive particle equations of motion from them
\end{enumerate}
These two approaches do not necessarily need to lead to the same particle equations of motion. Point 2. will definitely give us a structure-preserving equations of motion, but these might not be the "physically correct" trajectories, \textit{a priori}. For the first point, we need to understand how to get the particle equations of motion from \eqref{final-curvy-vlasov}, which we will discuss in the next section. For the second point, we need to figure out how to discretize the Hamiltonian structure, in particular, how to discretize the Poisson bracket, which is done in sec. \eqref{sec:discrete-hamiltonian-structure}. A more general discussion on the topics and methods of this section can be found in the documentation of \texttt{STRUPHY} \cite{StruphyDoc}. From now on we will only be dealing with equations on the logical domain and so we will suppress the hats to improve readability.

\subsection{Discrete Particle Equations of Motion}

We consider a transport equation with a source term in curvilinear coordinates
\begin{equation}\label{transport-equation}
	\del_t f_1 + \bA_1(\bh, \bv) \cdot \nabla_\bh f_1 + \bA_2(\bh, \bv) \cdot \nabla_\bv f_1 = S(\bh, \bv) \, ,
\end{equation}
where we assume that always
\begin{align}
	\nabla_\bh \cdot \left(\sqrt{g} \bA_1(\bh,\bv)\right) & = 0 \, , & \nabla_\bv \cdot \bA_2(\bh,\bv) & = 0 \, ,
\end{align}
as is the case for \eqref{final-cart-vlasov-equation} and \eqref{final-curvy-vlasov}.

Computing integrals of arbitrary functions $h$ against the distribution function $f_1$ can be done either by directly plugging in the PIC approximation
\begin{equation}\label{PIC-for-distribution-function}
	f_1(t, \bh, \bv) \approx f_{1,h}(t, \bh, \bv) = \frac{1}{N} \sum_p w_p(t) \; \delta(\bh - \bh_p(t)) \; \delta(\bv - \bv_p(t)) \, ,
\end{equation}
or multiplying the integrand by $\frac{s(\bh_p(t),\bv_p(t))}{s(\bh_p(t),\bv_p(t))}$, with $s$ a time-independent probability distribution function that is constant along the characteristics, such that the weights are related to $s$ by
\begin{equation}
	w_p(t) \equiv \frac{f_1(t,\bh_p(t),\bv_p(t))}{s(t,\bh_p(t),\bv_p(t))} = \frac{f_1(t,\bh_p(t),\bv_p(t))}{s(t=0,\bh_p(t=0),\bv_p(t=0))} = \frac{f_1(t,\bh_p(t),\bv_p(t))}{s_{0,p}}
\end{equation}
and are time-dependent because $f_1$ is not conserved along the characteristics. In both cases, the particles' positions $\bh_p$ and velocities $\bv_p$ are sampled from the coefficient of the 3-form $\sqrt{g} \, s(\bh, \bv)$.

The equations of motion for the particles are derived by demanding that \eqref{PIC-for-distribution-function} shall also verify \eqref{transport-equation} in a weak sense. We start by looking at the time derivative of the integral of $f_1(t,\bh,\bv)$ multiplied by a test function $h(\bh,\bv)$:
\begin{subequations}
	\begin{align}
		& \dt \int f_1(t,\bh,\bv) h(\bh,\bv) \sqrt{g} \d \bh\dv = \int h(\bh,\bv) \del_t f_1(t, \bh, \bv) \sqrt{g} \d \bh\dv \\
		& \quad = \int h(\bh, \bv) \left[ - \bA_1 \cdot \nabla_\bh f_1(t,\bh,\bv) - \bA_2 \cdot \nabla_\bv f_1(t,\bh,\bv) + S(\bh,\bv) \right] \sqrt{g} \d \bh \dv \\
		& \quad \approx \int f_{1, h} (t, \bh,\bv) \left[ \bA_1 \cdot \nabla_\bh h(\bh,\bv) + \bA_2 \cdot \nabla_\bv h(\bh,\bv) \right] \sqrt{g} \, \frac{s(t, \bh, \bv)}{s(t, \bh, \bv)} \d \bh \dv \\
		& \quad \quad + \int S(\bh, \bv) \, \frac{h(\bh,\bv)}{s(t, \bh, \bv)}{s(t, \bh,\bv)} \sqrt{g} \d \bh \dv \label{approximate-by-monte-carlo} \\
		& \quad = \frac{w_p(t)}{N} \left[ \bA_1(\bh_p(t), \bv_p(t)) \cdot \nabla_\bh h(\bh_p(t), \bv_p(t)) + \bA_2(\bh_p(t), \bv_p(t)) \cdot \nabla_\bv h(\bh_p(t), \bv_p(t)) \right] \\
		& \quad \quad + \frac{1}{N} \sum_p S(\bh_p(t), \bv_p(t)) \, \underbrace{\frac{h(\bh_p(t), \bv_p(t))}{s(t, \bh_p(t), \bv_p(t))}}_{=\frac{h(\bh_p(t), \bv_p(t))}{s_{0,p}}} \, ,
	\end{align}
\end{subequations}
where we used the Monte Carlo method to approximate the integral in line \eqref{approximate-by-monte-carlo}. Comparing this with the expression obtained by directly plugging in \eqref{PIC-for-distribution-function}
\begin{subequations}
	\begin{align}
		& \dt \int f_{1,h}(t,\bh,\bv) h(\bh,\bv) \sqrt{g} \d \bh\dv = \dt \frac{1}{N} \sum_p w_p(t) h(\bh_p(t), \bv_p(t)) \\
		& = \frac{1}{N} \sum_p \left[ h(\bh_p(t), \bv_p(t)) \diff{w_p}{t} + w_p(t) \nabla_\bh h(\bh_p(t),\bv_p(t)) \cdot \diff{\bh_p}{t} + w_p(t) \nabla_\bv h(\bh_p(t),\bv_p(t)) \cdot \diff{\bv_p}{t} \right] \, ,
	\end{align}
\end{subequations}
we can read off the equations of motion for the particle positions $\bh_p$, velocities $\bv_p$, and weights $w_p$
\begin{subequations}
	\begin{align}
		\diff{\bh_p}{t} & = \bA_1(\bh_p(t), \bv_p(t)) \, , \\
		\diff{\bv_p}{t} & = \bA_2(\bh_p(t), \bv_p(t)) \, , \\
		\diff{w_p}{t} & = \frac{1}{s_{0,p}} S(\bh_p(t), \bv_p(t)) \, .
	\end{align}
\end{subequations}
Applying this to \eqref{final-curvy-vlasov}, the equations of motion for the particle positions, velocities, and weights read:
\begin{subequations}\label{particle-eom}
	\begin{align}
		\diff{\bh_p}{t} & = DL^{-1} \bv_p(t) \, , \\
		\diff{\bv_p}{t} & = \frac{1}{\epsilon} DL^{-T} \bE_0(\bh_p(t)) + \frac{1}{\epsilon} \bv_p(t) \times \frac{1}{\sqrt{g}} \, DL \, \bB_0(\bh_p(t)) \, , \\
		\diff{w_p}{t} & = \frac{1}{s_{0,p}} \frac{1}{\vth^2 \, \epsilon} \left[DL^{-T} \bE_1(\bh_p(t))\right] \cdot \bv_p(t) f_{0,p} \, ,
	\end{align}
\end{subequations}
where we abbreviated $f_{0,p} = f_0(\bh_p(t), \bv_p(t))$ and also $s_{0,p} = s(\bh_p(t=0), \bv_p(t=0))$.

\subsection{Discrete Hamiltonian Structure}\label{sec:discrete-hamiltonian-structure}

We now discretize the Hamiltonian framework given by \eqref{final-curvy-Hamiltonian} and \eqref{final-curvy-bracket}.

\subsubsection{Discrete Conserved Hamiltonian}\label{sec:semi-discrete-hamiltonian}

We recall the Hamiltonian \eqref{final-curvy-Hamiltonian}
\begin{equation}
	H = \frac{\alpha^2}{2} \vth^2 \int_{\hat{\Omega}\times \real^3} \frac{f_1^2}{f_0} \sqrt{g}\d\bh\d\bv + H_\text{Maxwell} \, .
\end{equation}
The first term can be discretized using the Monte Carlo method
\begin{subequations}
	\begin{align}
		\int \frac{f_1^2}{f_0} \sqrt{g} \d\bh \dv & = \int \frac{\left(f_1(t, \bh, \bv)\right)^2}{f_0(\bh, \bv)} \frac{s(t, \bh, \bv)}{s(t, \bh,\bv)} \sqrt{g} \d\bh \dv \\
		& \approx \frac{1}{N} \sum_p \frac{\left(f_1(t,\bh_p, \bv_p)\right)^2}{s(t, \bh_p, \bv_p) \, f_0(\bh_p, \bv_p)} \\
		& = \frac{1}{N} \sum_p \frac{s_{0,p} \, w_p(t)^2}{f_{0,p}} \, .
	\end{align}
\end{subequations}
We will discuss the discretization of the fields in more detail below, so for now we replace $\bE_1 \to \bE_h$ and $\bB_1 \to \bB_h$. We note that the weak Ampère equation \eqref{curvy-Ampere} is discretized as:
\begin{subequations}\label{ad-hoc-weak-Ampere}
\begin{align}
		\dt \int \bF_h \cdot G^{-1} \bE_h \sqrt{g}\d\bh & = \int \left(\nabla_\bh\times \bF_h\right) \cdot G \, \bB_h \frac{1}{\sqrt{g}} \d\bh - \frac{\alpha^2}{\epsilon} \int DL^{-T} \bF_h\cdot \bv f_{1, h} \, \sqrt{g} \d\bh \dv \\
		& = \int \left(\nabla_\bh\times \bF_h\right) \cdot G \, \bB_h \frac{1}{\sqrt{g}} \d\bh - \frac{\alpha^2}{N \, \epsilon} \sum_p DL^{-T}(\bh_p) \bF_h(\bh_p) \cdot \bv_p \, w_p \, ,
\end{align}
\end{subequations}
where we have plugged in the PIC approach \eqref{PIC-for-distribution-function} for the distribution function $h$.

\subsubsection{Discrete Poisson Bracket}

In order to discretize the bracket \eqref{final-curvy-bracket}  we need to generalize the discussion in \cite{GEMPIC} in order to figure out what the discrete versions of
\begin{align}
	\var{\cF}{f_1} \, , && \nabla_\bh \var{\cF}{f_1} \, , && \text{and} && \nabla_\bv \var{\cF}{f_1}
\end{align}
are. We demand that the variation of a functional $\cF[f_1]$ shall not change under the discretization (at least in a weak sense). We also know that such a functional will become a function $F[w_p, \bh_p, \bv_p]$ of the particle positions, velocities, and weights. We demand that the variation of both these quantities be the same
\begin{equation}\label{variational-demand}
	\delta \cF[f_1] = \delta F[w_p, \bh_p, \bv_p] \, .
\end{equation}
The left-hand side reads:
\begin{subequations}
	\begin{align}
		\delta \cF[f_1] & = \int \var{\cF}{f_1} \delta f_{1,h} \sqrt{g} \d\bh \dv \\
		& = \int \var{\cF}{f_1} \sum_p \left[ \pa{f_{1,h}(w_p, \bh_p, \bv_p)}{w_p} \delta w_p + \nabla_{\bh_p} f_{1,h}(w_p, \bh_p, \bv_p) \cdot \delta \bh_p \right. \\
		& \hspace{22mm} \left. + \nabla_{\bv_p} f_{1,h}(w_p, \bh_p, \bv_p) \cdot \delta \bv_p \vphantom{\pa{h_h}{w_p}}\right] \sqrt{g} \d\bh \dv \\
		& = \int \var{\cF}{f_1} \frac{1}{N} \sum_p \left[ \delta(\bh-\bh_p) \delta(\bv-\bv_p) \delta w_p + w_p \delta(\bv-\bv_p) \nabla_{\bh_p} \delta(\bh - \bh_p) \cdot \delta \bh_p \right. \\
		& \hspace{25mm} \left. + w_p \delta(\bh-\bh_p) \nabla_{\bv_p} \delta(\bv-\bv_p) \cdot \delta \bv_p \right] \sqrt{g} \d\bh \dv \\
		& = \int \var{\cF}{f_1} \frac{1}{N} \sum_p \left[ \delta(\bh-\bh_p) \delta(\bv-\bv_p) \delta w_p - w_p \delta(\bv-\bv_p) \nabla_\bh \delta(\bh - \bh_p) \cdot \delta \bh_p \right. \\
		& \hspace{25mm} \left. - w_p \delta(\bh-\bh_p) \nabla_\bv \delta(\bv-\bv_p) \cdot \delta \bv_p \right] \sqrt{g} \d\bh \dv \\
		& = \frac{1}{N} \sum_p \int \left[ \var{\cF}{f_1} \delta w_p + w_p \nabla_\bh\var{\cF}{f_1} \cdot \delta \bh_p + w_p \nabla_\bv \var{\cF}{f_1} \cdot \delta \bv_p \right] \delta(\bh - \bh_p) \delta(\bv - \bv_p) \sqrt{g} \d\bh \dv \\
		& = \frac{1}{N} \sum_p \left[ \var{\cF}{f_1}(\bh_p,\bv_p) \delta w_p + w_p \nabla_\bh\var{\cF}{f_1}(\bh_p,\bv_p) \cdot \delta \bh_p + w_p \nabla_\bv \var{\cF}{f_1}(\bh_p,\bv_p) \cdot \delta \bv_p \right] \, ,
	\end{align}
\end{subequations}
where we used that
\begin{equation}
	\nabla_\bh \delta(\bh - \bh_p) = - \nabla_{\bh_p} \delta(\bh - \bh_p) \, .
\end{equation}
The right-hand side of \eqref{variational-demand} reads
\begin{equation}
		\delta F[w_p, \bh_p, \bv_p] = \sum_p \pa{F}{w_p} \delta w_p + \pa{F}{\bh_p} \cdot \delta \bh_p + \pa{F}{\bv_p} \cdot \delta \bv_p \, ,
\end{equation}
and by comparing terms, we find the following replacement rules
\begin{align}\label{replacement-rules}
	\frac{1}{N} \var{\cF}{f_1}(\bh_p, \bv_p) & = \pa{F}{w_p} \, , & \frac{w_p}{N} \nabla_\bh \var{\cF}{f_1}(\bh_p, \bv_p) & = \pa{F}{\bh_p} \, , & \frac{w_p}{N} \nabla_\bv \var{\cF}{f_1}(\bh_p, \bv_p) & = \pa{F}{\bv_p} \, .
\end{align}
Thus, we can discretize \eqref{final-curvy-bracket} using the Monte Carlo method, and plugging in these replacement rules yields
\begin{equation}\label{discrete-Poisson-bracket}
	\begin{aligned}
		\bracket{F}{G} & = \int \left[ \left( \nabla_\bh \times \var{F}{\bE_h} \right) \cdot G \, \var{G}{\bB_h} - \left( \nabla_\bh \times \var{G}{\bE_h} \right) \cdot G \, \var{F}{\bE_h} \right] \frac{1}{\sqrt{g}} \d\bh \\
		& \quad + \frac{1}{\vth^2 \, \epsilon} \sum_p \frac{f_{0,p}}{s_{0,p}} \left[ \pa{F}{w_p} \var{G}{\bE_h}(\bh_p, \bv_p) - \pa{G}{w_p} \var{F}{\bE_h} (\bh_p, \bv_p) \right] \cdot DL^{-1} \bv_p \\
		& \quad \begin{aligned}
			& + \frac{N}{\alpha^2 \vth^2} \sum_p \frac{f_{0,p}}{s_{0,p} \, w_p} \left[ DL^{-1} \bv_p \cdot \left(\pa{F}{\bh_p} \pa{G}{w_p} - \pa{F}{w_p} \pa{G}{\bh_p} \right) \right. \\
			& \hspace{35mm} \left. + \frac{1}{\epsilon} \left(DL^{-T} \bE_0 + \bv_p \times \sqrt{g} DL\, \bB_0\right) \cdot \left(\pa{F}{\bv_p} \pa{G}{w_p} - \pa{F}{w_p} \pa{G}{\bv_p} \right) \right] \, .
		\end{aligned}
	\end{aligned}
\end{equation}
We note that in comparison to \eqref{final-curvy-bracket}, we have removed the factor of $\frac{1}{2}$ in the last term, which was there to produce the correct equations of motion. Since the notion of integration by parts does not exist in the discrete setting, this factor would still appear in the equations of motion, and hence we removed it. It can be easily checked that this bracket and the Hamiltonian from sec. \eqref{sec:semi-discrete-hamiltonian} produces the correct particle equations of motion.

% !TeX spellcheck = en_GB

\section{Field Discretization}\label{sec:feec-discretization}

We will follow \cite{GEMPIC, Holderied2021} in discretizing the electromagnetic fields using finite element exterior calculus (FEEC). This section mainly introduces notation; for details, we refer the reader to the mentioned works.

\subsection{FEEC}

The heart of the FEEC method is the following diagram
\begin{equation}
	\begin{tikzcd}[row sep=large]
		H^1(\hat{\Omega}) \arrow[r,"\;\hnab_\bh\;"] \arrow[d, "\Pi_0" left] & H(\tcurl, \hat{\Omega}) \arrow[r,"\;\hnab_\bh\times\;"] \arrow[d, "\Pi_1" left] & H(\tdiv,\hat{\Omega}) \arrow[r,"\;\hnab_\bh\cdot\;"] \arrow[d, "\Pi_2" left]  & L^2(\hat{\Omega}) \arrow[d, "\Pi_3" left] \\
		V^0_h \arrow[r, "\hnab_\bh"] & V^1_h \arrow[r, "\hnab_\bh \times"] & V^2_h \arrow[r, "\hnab_\bh\cdot"] & V^3_h
	\end{tikzcd} \, ,
\end{equation}
which is comprised of the de Rham sequence \eqref{logical-de-Rham-sequence} on the logical domain and a corresponding de Rham sequence on the finite element (FE) approximation spaces $V_i$. These two are linked by the projection operators $\Pi_k$.

We follow the notation of \cite{GEMPIC, Holderied2021} by denoting the basis functions using $\bL^k$ for scalar-valued spaces and $\fL^k$ for vector-valued spaces, where $k$ indicates in which FE space $V_h^k$ we are. So the charge density $\rho$ lives in $H^1(\hat{\Omega})$; its FE approximation will be written as
\begin{equation}
	H^1(\hat{\Omega}) \ni \quad \rho(\bh) \approx \rho_h(\bh) = \br^T \bL^0(\bh) \quad \in V^0_h \, ,
\end{equation}
and analogously for any function $f$ that lives in $L^2(\hat{\Omega})$
\begin{equation}
	L^2(\hat{\Omega}) \ni \quad f(\bh) \approx f_h(\bh) = \bf^T \bL^3(\bh) \quad \in V^3_h \, .
\end{equation}
For the vector-valued spaces, we have, for example, the electric field $\bE \in H(\tcurl, \hat{\Omega})$ and $\bB \in H(\tdiv, \hat{\Omega})$, which are approximated as follows
\begin{subequations}
	\begin{align}
		H(\tcurl, \hat{\Omega}) \ni \quad \bE_1(\bh) \approx \left( \bE_h(\bh) \right)^T & = \be_1^T \fL^1(\bh) \quad \in V^1_h \, , \\
		H(\tdiv, \hat{\Omega}) \ni \quad \bB_1(\bh) \approx \left( \bB_h(\bh) \right)^T & = \bb_1^T \fL^2(\bh) \quad \in V^2_h \, .
	\end{align}
\end{subequations}
We also denote mass matrices in these spaces by $\fM_k$.

\subsection{The Full Semi-discrete System}

After the FEEC discretization, the particle equations of motion \eqref{particle-eom} become
\begin{subequations}\label{discrete-eoms}
	\begin{align}
		\dt \bh_p & = DL^{-1}(\bh_p) \bv_p \, , \\
		\dt \bv_p & = \frac{1}{\epsilon} DL^{-T}(\bh_p) (\fL^1(\bh_p))^T \be_0 + \frac{1}{\epsilon} \bv_p \times \frac{1}{\sqrt{g(\bh_p)}} \, DL(\bh_p) (\fL^2(\bh_p))^T \bb_0 \, , \\
		\dt w_p & = \frac{1}{s_0} \frac{1}{\vth^2 \epsilon} \left[DL^{-T}(\bh_p) (\fL^1(\bh_p))^T \be_1\right] \cdot \bv_p \, f_{0,p} \, .
	\end{align}
\end{subequations}
The Ampère equation (in weak form) \eqref{ad-hoc-weak-Ampere} becomes
\begin{equation}\label{discrete-ampere}
	\dt \fM_1 \be_1 = \fC^T \fM_2 \bb_1 - \frac{\alpha^2}{N \, \epsilon} \sum_p w_p \fL^1(\bh_p) DL^{-1} (\bh_p) \bv_p \, ,
\end{equation}
and Faraday's law (in strong form) \eqref{curvy-Faraday} becomes
\begin{equation}\label{discrete-faraday} \, .
	\dt \bb_1 = - \fC \be_1
\end{equation}
The Hamiltonian \eqref{final-curvy-Hamiltonian} is finally discretized as
\begin{equation}\label{final-discrete-Hamiltonian}
	H = \frac{1}{2N} \alpha^2 \vth^2 \sum_p \frac{s_{0,p} \, w_p^2}{f_{0,p}} + \frac{1}{2} \be_1^T \fM_1 \be_1 + \frac{1}{2} \bb_1^T \fM_2 \bb_1 \, ,
\end{equation}
and the Poisson bracket \eqref{discrete-Poisson-bracket} becomes (under the replacement rules found in \cite{GEMPIC})
\begin{subequations}\label{final-discrete-Poisson-bracket}
	\begin{align}
		\bracket{F}{G}
		& = \pa{F}{\be_1} \fM_1^{-1} \fC^T \pa{G}{\bb_1} - \pa{F}{\bb_1} \fC \fM_1^{-1} \pa{G}{\be_1} \\
		& \quad + \frac{1}{\vth^2 \, \epsilon} \sum_p \frac{f_{0,p}}{s_{0,p}} \left[ \pa{F}{w_p} DL^{-1} \bv_p \cdot (\fL^1)^T \fM_1^{-1} \pa{G}{\be_1} - \pa{F}{\be_1} \fM^{-1} \fL^1 \cdot DL^{-1} \bv_p \pa{G}{w_p} \right] \\
		& \quad + \frac{N}{\alpha^2 \vth^2} \sum_p \frac{f_{0,p}}{s_{0,p} w_p} \left[ DL^{-1} \bv_p \cdot \left(\pa{F}{\bh_p} \pa{G}{w_p} - \pa{F}{w_p} \pa{G}{\bh_p} \right) \right. \\
		& \hspace{20mm} \left. + \frac{1}{\epsilon} \left(DL^{-T} (\fL^1)^T \be_0 + \bv_p \times \frac{1}{\sqrt{g}} DL\, (\fL^2)^T \bb_0\right) \cdot \left(\pa{F}{\bv_p} \pa{G}{w_p} - \pa{F}{w_p} \pa{G}{\bv_p} \right) \right] \, ,
	\end{align}
\end{subequations}
where we have suppressed all arguments for increased readability.

The bracket \eqref{final-discrete-Poisson-bracket} is in fact a Poisson bracket as it is manifestly anti-symmetric and verifies the Leibniz rule. It remains to show the Jacobi identity
\begin{equation}\label{LVM-jacobi}
	\bracket{\bracket{U_i}{U_j}}{U_k} + \bracket{\bracket{U_k}{U_i}}{U_j} + \bracket{\bracket{U_j}{U_k}}{U_i} = 0 \, ,
\end{equation}
for $U_i \in \left\{ (\bh_p)_i, (\bv_p)_i, w_p, (\be_1)_i, (\bb_1)_i \right\}$. We only consider non-vanishing, non-constant brackets. We also note that
\begin{equation}
	\bracket{U_p}{U_q} = \delta_{p,q} \qquad \qquad \forall \, U \in \left\{w, \bh, \bv\right\} \, ,
\end{equation}
i.e., the particle part of the bracket vanishes if the arguments are not from the same particle. We start with
\begin{equation}
	\bracket{w_p}{\be_1} = - \bracket{\be_1}{w_p} = \frac{1}{\vth^2 \epsilon} \frac{f_{0,p}}{s_{0,p}} \, \bv_p \cdot \left( DL^{-T} \fL^1 \fM_1^{-1} \right)
\end{equation}
which is a function of $\bh_p$ and $\bv_p$, thus only a further bracket with $w_p$ does not vanish. Then the Jacobi identity reads
\begin{equation}
	\bracket{\bracket{w_p}{\be_1}}{w_p} + \bracket{\bracket{w_p}{w_p}}{\be_1} + \bracket{\bracket{\be_1}{w_p}}{w_p} = \bracket{\bracket{w_p}{\be_1}}{w_p} - \bracket{\bracket{w_p}{\be_1}}{w_p} = 0
\end{equation}
because of the antisymmetry of the bracket. Next, we consider
\begin{equation}
	\bracket{w_p}{\bh_p} = - \bracket{\bh_p}{w_p} = - \frac{N}{\alpha^2 \vth^2} \frac{f_{0,p}}{s_{0,p} w_p} \, DL^{-1} \bv_p
\end{equation}
which vanishes by the same argument. Finally
\begin{equation}
	\bracket{w_p}{\bh_p} = - \bracket{\bh_p}{w_p} = - \frac{N}{\alpha^2 \vth^2} \frac{f_{0,p}}{s_{0,p} w_p} \frac{1}{\epsilon} \left(DL^{-T} \fL^1 \be_0 + \bv_p \times \left(\frac{1}{\sqrt{g}} DL\, \fL^2 \bb_0\right) \right)
\end{equation}
and also this gives a zero Jacobi identity, so \eqref{LVM-jacobi} is verified for all coordinate functions and therefore for arbitrary functions of them.

We show that the bracket \eqref{final-discrete-Poisson-bracket} is in fact a Poisson bracket. Antisymmetry and Leibniz rule are again manifest from the construction. We still need to show the Jacobi identity
\begin{equation}\label{LVM-jacobi2}
	\bracket{\bracket{U_i}{U_j}}{U_k} + \bracket{\bracket{U_k}{U_i}}{U_j} + \bracket{\bracket{U_j}{U_k}}{U_i} = 0 \, ,
\end{equation}
for $U_i \in \left\{ (\bh_p)_i, (\bv_p)_i, w_p, (\be_1)_i, (\bb_1)_i \right\}$. We only consider non-vanishing, non-constant brackets. We also note that
\begin{equation}
	\bracket{U_p}{U_q} = \delta_{p,q} \qquad \qquad \forall \, U \in \left\{w, \bh, \bv\right\} \, ,
\end{equation}
i.e., the particle part of the bracket vanishes if the arguments are not from the same particle.

We start with
\begin{equation}
	\bracket{w_p}{\be_1} = - \bracket{\be_1}{w_p} = \frac{1}{\vth^2 \epsilon} \frac{f_{0,p}}{s_{0,p}} \, \bv_p \cdot \left( DL^{-T} \fL^1 \fM_1^{-1} \right) \, ,
\end{equation}
which is a function of $\bh_p$ and $\bv_p$, thus only a further bracket with $w_p$ does not vanish. Then the Jacobi identity reads
\begin{equation}
	\bracket{\bracket{w_p}{\be_1}}{w_p} + \bracket{\bracket{w_p}{w_p}}{\be_1} + \bracket{\bracket{\be_1}{w_p}}{w_p} = \bracket{\bracket{w_p}{\be_1}}{w_p} - \bracket{\bracket{w_p}{\be_1}}{w_p} = 0 \, ,
\end{equation}
because of the antisymmetry of the bracket. Next, we consider
\begin{subequations}
	\begin{align}
		\bracket{w_p}{\bh_p} & = - \bracket{\bh_p}{w_p} = - \frac{N}{\alpha^2 \vth^2} \frac{f_{0,p}}{s_{0,p} w_p} \, DL^{-1} \bv_p \, , \\
		\bracket{w_p}{\bv_p} & = - \bracket{\bv_p}{w_p} = - \frac{N}{\alpha^2 \vth^2} \frac{f_{0,p}}{s_{0,p} w_p} \frac{1}{\epsilon} \left(DL^{-T} \fL^1 \be_0 + \bv_p \times \left(\frac{1}{\sqrt{g}} DL\, \fL^2 \bb_0\right) \right) \, ,
	\end{align}
\end{subequations}
which depend on all $\bh_p$, $\bv_p$, and $w_p$, and thus have non-vanishing brackets with $\bh_p$ and $\bv_p$. The following Jacobi identity terms
\begin{subequations}
	\begin{align}
		\bracket{\bracket{w_p}{\bh_p}}{w_p} + \bracket{\bracket{w_p}{w_p}}{\bh_p} + \bracket{\bracket{\bh_p}{w_p}}{w_p} = \bracket{\bracket{w_p}{\bh_p}}{w_p} - \bracket{\bracket{w_p}{\bh_p}}{w_p} & = 0 \, , \\
		\bracket{\bracket{w_p}{\bh_p}}{\bh_p} + \bracket{\bracket{\bh_p}{w_p}}{\bh_p} + \bracket{\bracket{\bh_p}{\bh_p}}{w_p} = \bracket{\bracket{w_p}{\bh_p}}{\bh_p} - \bracket{\bracket{w_p}{\bh_p}}{\bh_p} & = 0 \, ,
	\end{align}
\end{subequations}
as well as
\begin{subequations}
	\begin{align}
		\bracket{\bracket{w_p}{\bv_p}}{w_p} + \bracket{\bracket{w_p}{w_p}}{\bv_p} + \bracket{\bracket{\bv_p}{w_p}}{w_p} = \bracket{\bracket{w_p}{\bv_p}}{w_p} - \bracket{\bracket{w_p}{\bv_p}}{w_p} & = 0 \, , \\
		\bracket{\bracket{w_p}{\bv_p}}{\bv_p} + \bracket{\bracket{\bv_p}{w_p}}{\bv_p} + \bracket{\bracket{\bv_p}{\bv_p}}{w_p} = \bracket{\bracket{w_p}{\bv_p}}{\bv_p} - \bracket{\bracket{w_p}{\bv_p}}{\bv_p} & = 0 \, ,
	\end{align}
\end{subequations}
vanish by the same argument as above. However, there are also mixed terms
\begin{equation}
	\bracket{\bracket{w_p}{\bh_p}}{\bv_p} + \bracket{\bracket{\bv_p}{w_p}}{\bh_p} + \bracket{\bracket{\bh_p}{\bv_p}}{w_p} \, ,
\end{equation}
which could be non-vanishing. To check this, we have to compute the nested brackets
\begin{subequations}
	\begin{align}
		\bracket{\bracket{w_p}{\bh_p}}{\bv_p} & = - \frac{N}{\alpha^2 \vth^2} \frac{f_{0,p}}{s_{0,p} w_p} \frac{1}{\epsilon} \left(DL^{-T} \fL^1 \be_0 + \bv_p \times \left(\frac{1}{\sqrt{g}} DL\, \fL^2 \bb_0\right) \right) \cdot \pa{}{w_p} \bracket{w_p}{\bh_p} \\
		& = - \frac{N}{\alpha^2 \vth^2} \frac{f_{0,p}}{s_{0,p} w_p} \frac{1}{\epsilon} \left(DL^{-T} \fL^1 \be_0 + \bv_p \times \left(\frac{1}{\sqrt{g}} DL\, \fL^2 \bb_0\right) \right) \cdot \left(\frac{N}{\alpha^2 \vth^2} \frac{f_{0,p}}{s_{0,p} w_p^2} \, DL^{-1} \bv_p\right) \, , \\
		\bracket{\bracket{\bv_p}{w_p}}{\bh_p} & = \frac{N}{\alpha^2 \vth^2} \frac{f_{0,p}}{s_{0,p} w_p} \, \left(DL^{-1} \bv_p\right) \cdot \pa{}{w_p} \bracket{\bv_p}{w_p} \\
		& = \frac{N}{\alpha^2 \vth^2} \frac{f_{0,p}}{s_{0,p} w_p} \, \left(DL^{-1} \bv_p\right) \cdot \left(\frac{N}{\alpha^2 \vth^2} \frac{f_{0,p}}{s_{0,p} w_p^2} \frac{1}{\epsilon} \left(DL^{-T} \fL^1 \be_0 + \bv_p \times \left(\frac{1}{\sqrt{g}} DL\, \fL^2 \bb_0\right) \right)\right) \, ,
	\end{align}
\end{subequations}
and so we see immediately that
\begin{equation}
	\bracket{\bracket{w_p}{\bh_p}}{\bv_p} = - \bracket{\bracket{\bv_p}{w_p}}{\bh_p} \, ,
\end{equation}
and so
\begin{equation}
	\bracket{\bracket{w_p}{\bh_p}}{\bv_p} + \bracket{\bracket{\bv_p}{w_p}}{\bh_p} + \bracket{\bracket{\bh_p}{\bv_p}}{w_p} = 0 \, ,
\end{equation}
because $\bracket{\bh_p}{\bv_p} = 0$. Therefore \eqref{LVM-jacobi2} is verified for all coordinate functions and thus for arbitrary functions of them.

% !TeX spellcheck = en_GB

\section{Poisson Splitting \& Time Discretization}\label{sec:poisson-splitting}

We employ a Poisson splitting of the matrix \eqref{final-discrete-Poisson-bracket} in order to preserve the Hamiltonian \eqref{final-discrete-Hamiltonian} exactly in each substep. We split five terms and subsequently arrive at five substeps, which we can combine to achieve any order in time. The framework of \texttt{STRUPHY} \cite{StruphyPaper, StruphyDoc} allows us to choose flexibly between the Lie-Trotter splitting (a first-order method in time) or Strang splitting (a second-order method).

\subsection{Maxwell Substep}

The first term in the bracket \eqref{final-discrete-Poisson-bracket} reads
\begin{equation}
	\bracket{\cF}{\cG}_1 = \pa{\cF}{\be_1} \fM_1^{-1} \fC^T \pa{\cG}{\bb_1} - \pa{\cF}{\bb_1} \fC \fM_1^{-1} \pa{\cG}{\be_1} \, ,
\end{equation}
and yields the discrete Maxwell equations in vacuum
\begin{subequations}
	\begin{align}
		\dt \fM_1 \be_1 & = \fC^T \fM_2 \bb_1 \, , \\
		\dt \bb_1 & = - \fC \be_1 \, .
	\end{align}
\end{subequations}
This system of equations can be written as a $2\times 2$ system
\begin{subequations}
	\begin{align}
		\dt \begin{pmatrix} \fM_1 \be_1 \\ \bb_1 \end{pmatrix} =
		\begin{pmatrix}
			0 & \fC^T \fM_2 \\
			- \fC & 0
		\end{pmatrix} \begin{pmatrix} \be_1 \\ \bb_1 \end{pmatrix} \, ,
	\end{align}
\end{subequations}
which can be solved using a Crank-Nicolson scheme as the discretization in time and the Schur method \cite{Haynsworth_1968}. Since the Hamiltonian is at most quadratic in the variables, the Crank-Nicolson scheme is a discrete gradient and therefore energy-preserving.

The update for $\be_1$ then reads
\begin{equation}\label{e-update}
	(M/D) \be_1^{n+1} = \left( \fM_1 + \fA\fH \right) \be_1^n - 2\fA \bb_1^n \, ,
\end{equation}
with the matrices
\begin{align}
	\fA & = - \frac{\Dt}{2} \fC^T \fM_2 \, , \\
	\fA \fH & = - \frac{\Dt^2}{4} \fC^T \fM_2 \fC \, , \\
	(M/D) & = \fM_1 - \fA\fH = \fM_1 + \frac{\Dt^2}{4} \fC^T \fM_2 \fC \, ,
\end{align}
and the update for $\bb_1$ reads
\begin{equation}
	\bb_1^{n+1} = \bb_1^n - \frac{\Dt}{2} \fC \left( \be_1^{n+1} + \be_1^n \right) \, .
\end{equation}
Note that for solving \eqref{e-update}, we have to invert the matrix $(M/D)$, and thus energy conservation can only be guaranteed up to the accuracy of the inversion solver.

\subsection{Coupling Substep}\label{sec:subsystem-two}

The second term of the Poisson bracket,
\begin{equation}
	\bracket{\cF}{\cG}_2 = \frac{1}{\vth^2 \, \epsilon} \sum_p \frac{f_{0,p}}{s_{0,p}} \left[ \pa{\cF}{w_p} DL^{-1} \bv_p \cdot (\fL^1)^T \fM_1^{-1} \pa{\cG}{\be_1} - \pa{\cF}{\be_1} \fM^{-1} \fL^1 \cdot DL^{-1} \bv_p \pa{\cG}{w_p} \right] \, ,
\end{equation}
couples the particles and the fields through the weights and the electric field. The equations read
\begin{subequations}\label{second-substep-initial-equations}
	\begin{align}
		\dt w_p & = \frac{f_{0,p}}{s_{0,p}} \frac{1}{\vth^2 \, \epsilon} \left[DL^{-T}(\bh_p) (\fL^1(\bh_p))^T \be_1\right] \cdot \bv_p \, , \\
		\dt \fM_1 \be_1 & = - \frac{\alpha^2}{N \, \epsilon} \sum_p w_p \, \fL^1(\bh_p) \cdot \left(DL^{-1} (\bh_p) \bv_p\right) \, .
	\end{align}
\end{subequations}
We can rewrite this by introducing a vector and matrices with sizes of the particle number:
\begin{align}
	\left(\bW\right)_p & \equiv w_p \, ,
	& \left(\fW\right)_p & \equiv \frac{f_{0,p}}{s_{0,p}} \frac{1}{\vth^2 \, \epsilon} \left(DL^{-1} \bv_p \right) \cdot \left(\fL^1\right)^T \, ,
	& \left(\fE \right)_p & \equiv \frac{\alpha^2}{N \, \epsilon} \fL^1 \cdot \left(DL^{-1} \bv_p \right) \, ,
\end{align}
which have sizes
\begin{align}
	\bW & \, \in \real^{N_p} \, , & \fW & \, \in \real^{N_p, N_1} \, , & \fE & \, \in \real^{N_1, N_p} \, .
\end{align}
\eqref{second-substep-initial-equations} can be written as a $2\times 2$ system
\begin{equation}
	\dt \begin{pmatrix} \fM_1 \be_1 \\ \bW \end{pmatrix} =
	\begin{pmatrix}
		0 & - \fE \\
		\fW & 0
	\end{pmatrix}
	\begin{pmatrix} \be_1 \\ \bW \end{pmatrix} \, ,
\end{equation}
so again, we can do an energy-conserving discretization in time using the Crank-Nicolson scheme and use the Schur complement \cite{Haynsworth_1968} to solve the resulting system. Here, the number of markers is usually much larger than the number of basis functions of $V_h^1$: $N_p \gg N_1$, so using the Schur complement shows its strength in computational cost while still being energy-conserving. We find that the update rule for the electric FE coefficients
\begin{equation}
	(M/D)\be_1^{n+1} = \left[\fM_1 - \frac{\Dt^2}{4} \fE \fW \right] \be_1^n - \Dt \fE \bW \, ,
\end{equation}
where
\begin{equation}
	(M/D) = \fM_1 - \frac{\Dt^2}{4} \fE \fW \, .
\end{equation}
The quantities $\fE \fW$ and $\fE \bW$ include summation over the particles, which is called accumulation in \texttt{STRUPHY} \cite{StruphyDoc}. Written out, they read
\begin{subequations}
	\begin{align}
		\fE \fW & = \frac{\alpha^2}{N \, \vth^2 \, \epsilon^2} \sum_p \frac{f_{0,p}}{s_{0,p}} \fL^1(\bh_p) \cdot \left(DL^{-1}(\bh_p) \bv_p \right) \left(DL^{-1}(\bh_p) \bv_p\right) \cdot \left(\fL^1(\bh_p)\right)^T \, , \\
		\fE \bW & = \frac{\alpha^2}{N \, \epsilon} \sum_p w_p \, \fL^1(\bh_p) DL^{-1}(\bh_p) \bv_p \, .
	\end{align}
\end{subequations}
Note that the accumulation matrix $\fE \fW$ is symmetric. Finally, $\bW^{n+1}$ is computed as
\begin{equation}
	\bW^{n+1} = \bW^n + \frac{\Dt}{2} \fW \left( \be^{n+1} + \be^n \right) \, ,
\end{equation}
which can be computed element-wise: for each component (i.e., each particle), we have
\begin{equation}
	w_p^{n+1} = w_p^n + \frac{\Dt}{2} \frac{f_{0,p}}{s_{0,p}} \frac{1}{\vth^2 \, \epsilon} \left( DL^{-1}(\bh_p) \bv_p \right) \cdot \left[ \left( \fL^1(\bh_p) \right)^T \left( \be_1^{n+1} + \be_1^n \right) \right] \, .
\end{equation}

\subsection{Space Advection Substep}\label{sec:subsystem-three}

The third term of the Poisson bracket reads
\begin{equation}
	\bracket{\cF}{\cG}_3 = \frac{N}{\alpha^2 \vth^2} \sum_p \frac{1}{s_{0,p} w_p} \left[ DL^{-1} \bv_p \cdot \left(\pa{\cF}{\bh_p} \pa{\cG}{w_p} - \pa{\cF}{w_p} \pa{\cG}{\bh_p} \right) \right] \, ,
\end{equation}
and yields the equations of motion
\begin{equation}
	\dt \bh_p = DL^{-1}(\bh_p) \bv_p \, .
\end{equation}
This can be solved very accurately using an explicit Runge-Kutta method of order 4, which will not conserve the energy exactly. For our use cases $n_0(\bx) = const.$, and therefore the Hamiltonian does not depend on the positions and is thus conserved during this step.

\subsection{Electric Lorentz Substep}

The fourth term of the bracket
\begin{equation}
	\bracket{\cF}{\cG} = \frac{N}{\alpha^2 \vth^2 \epsilon} \sum_p \frac{1}{s_{0,p} w_p} \left(DL^{-T} (\fL^1)^T \be_0 \right) \cdot \left(\pa{\cF}{\bv_p} \pa{\cG}{w_p} - \pa{\cF}{w_p} \pa{\cG}{\bv_p} \right) \, ,
\end{equation}
is the electric part of the Lorentz force of the background fields
\begin{equation}
	\dt \bv_p = \frac{1}{\epsilon} DL^{-T}(\bh_p) (\fL^1(\bh_p))^T \be_0 \, ,
\end{equation}
and can be solved analytically.

\subsection{Magnetic Lorentz Substep}

The fifth term of the Poisson bracket
\begin{equation}
	\bracket{\cF}{\cG} = \frac{N}{\alpha^2 \vth^2 \epsilon} \sum_p \frac{1}{s_{0,p} w_p} \left( \bv_p \times \frac{1}{\sqrt{g}} DL\, (\fL^2)^T \bb_0\right) \cdot \left(\pa{\cF}{\bv_p} \pa{\cG}{w_p} - \pa{\cF}{w_p} \pa{\cG}{\bv_p} \right) \, ,
\end{equation}
is the magnetic part of the Lorentz force of the background fields 
\begin{equation}
	\dt \bv_p = \frac{1}{\epsilon} \bv_p \times \frac{1}{\sqrt{g(\bh_p)}} \, DL(\bh_p) (\fL^2(\bh_p))^T \bb_0 \, .
\end{equation}
This system can also be solved analytically.

% !TeX spellcheck = en_GB

\section{Numerical Experiments}\label{sec:results}

We test our model by simulating the weak Landau damping and the dispersion relation under a constant magnetic background field. Our model exhibits superior properties for both cases when comparing it to a direct delta-f method. Both models are implemented in \texttt{STRUPHY} \cite{StruphyPaper, StruphyDoc} and the parameter files for the simulations are available in \cite{StruphySimulations}.

\subsection{Direct Delta-f Method}

The equations of motion are as in the text (eqns. \eqref{final-cart-maxwell-equations}), such that they can be discretized in the GEMPIC framework, yielding the same semi-discrete equations of motion as in the text \eqref{discrete-eoms}, \eqref{discrete-ampere}, and \eqref{discrete-faraday}. The source term is again responsible for dynamical weights, however, now they are split into a separate subsystems:
\begin{subequations}
	\begin{align}
		\text{(I): } & \left\{ \dt \bh_p = DL^{-1}(\bh_p) \bv_p \, , \right. \\
		\text{(II): } & \left\{ \dt \bv_p = \frac{1}{\epsilon} DL^{-T}(\bh_p) (\fL^1(\bh_p))^T \be_0 \, , \right. \\
		\text{(III): } & \left\{ \dt \fM_1 \be_1 = - \frac{\alpha^2}{N \, \epsilon} \sum_p w_p \fL^1(\bh_p) DL^{-1} (\bh_p) \bv_p \, , \right. \\
		\text{(IV): } & \left\{ \dt w_p = \frac{1}{\vth^2 \epsilon} \sum_p \bv_p \cdot (\fL^1(\bh_p))^T \be_1 \, , \right. \\
		\text{(V): } & \left\{ \dt \bv_p = \frac{1}{\epsilon} \bv_p \times \frac{1}{\sqrt{g(\bh_p)}} DL(\bh_p) (\fL^2(\bh_p))^T \left(\bb_0 + \bb_1\right) \, , \right. \\
		\text{(VI): } & \left\{ \begin{aligned}
			\dt \fM_1 \be_1 & = \fC^T \fM_2 \bb_1 \, , \\
			\dt \bb_1 & = - \fC \be_1 \, ,
		\end{aligned} \right.
	\end{align}
\end{subequations}
all of which can be solved either analytically or using the techniques discussed in sec. \eqref{sec:poisson-splitting}.

\subsection{Landau Damping}

Landau damping \cite{Landau_damping, Landau_1946} is a widely known effect in plasmas where an initial perturbation is damped exponentially by the particle-field interaction. The initial condition for the distribution functions reads
\begin{equation}
	f(t=0, \bx, \bv) = \frac{1}{\sqrt{2 \pi \vth^2}^3} \exp \left( - \frac{\abs{\bv}^2}{\vth^2} \right) \, \left[ 1 + \delta \cos\left(k x\right) \vphantom{\sum} \right] \, ,
\end{equation}
with parameters $\vth = 1$, $k=0.5$, and $\delta=10^{-3}$. The initial electric field is computed from the initial distribution function by solving a Poisson equation. The magnetic field is initialized as zero and should remain zero over time.

\subsubsection{Linear Vlasov--Ampère Model and Setup}

Since this is an electrostatic test case, we can reduce our model by setting $\bB_1 = 0$, hence we effectively solve a Vlasov--Ampère model. We also simulate this test case using the Vlasov--Maxwell model, which we discuss below. Furthermore, the background fields are identically zero, meaning the Vlasov equation only contains a transport in $\bh$. The units of scale are chosen such that $\alpha = 1$ and $\epsilon = -1$ (since it carries the sign of the charge; for the simulations, the sign does not play a role). So finally \eqref{final-curvy-vlasov} and \eqref{curvy-Ampere} become
\begin{subequations}
	\begin{align}
		\del_t \hat{f}_1 + \bv \cdot DL^{-T} \hnab_\bh \hat{f}_1 & = \frac{1}{\vth^2} DL^{-T} \hbE_1 \cdot \bv \hat{f}_0 \, , \\
		\dt \int_{\hat{\Omega}} \hbF \cdot G^{-1} \, \hbE_1 \, \sqrt{g} \d \bh & = \int_{\hat{\Omega}\times \real^3} DL^{-T} \, \hbF \cdot \bv \, \hat{f}_1 \, \sqrt{g} \d \bh \dv \, ,
	\end{align}
\end{subequations}
which leads to the following semi-discrete equations of motion
\begin{subequations}
	\begin{align}
		\dt \bh_p & = DL^{-1}(\bh_p) \bv_p \, , \\
		\dt \bv_p & = 0 \, , \\
		\dt w_p & = \frac{1}{s_{0,p}} \frac{1}{\vth^2} \left[DL^{-T}(\bh_p) (\fL^1(\bh_p))^T \be_1\right] \cdot \bv_p \, f_{0,p} \, , \\
		\dt \fM_1 \be & = \frac{1}{N} \sum_p w_p f_{0,p} \, \fL^1(\bh_p) DL^{-1} (\bh_p) \bv_p \, ,
	\end{align}
\end{subequations}
which we split and solve in time exactly as mentioned above. The remaining substeps are the coupling of weights and electric field \eqref{sec:subsystem-two} and the transport in position space by velocity \eqref{sec:subsystem-three}.

For our model, the initial conditions translate to
\begin{subequations}
	\begin{align}
		f_0(\bv) & = \frac{1}{\sqrt{2 \pi \vth^2}^3} \exp \left( - \frac{\abs{\bv}^2}{\vth^2} \right) \, , \\
		f_1(t=0, \bx, \bv) & = \frac{\delta}{\sqrt{2 \pi \vth^2}^3} \exp \left( - \frac{\abs{\bv}^2}{\vth^2} \right) \cos\left(k x\right) \, ,
	\end{align}
\end{subequations}
and the initial $\bE_1$ is computed from $f_1(t=0, \bx, \bv)$ using the Poisson equation.

Since we want to simulate the 1D case, the domain is chosen as $\Omega = [0, 4\pi] \times [0,1]^2$ and $k=\frac12$ such that the initial perturbation does exactly half a period on the domain. The numerical parameters are chosen as 32 grid points in $x$-direction and 1 in the other two spatial directions. The FE basis functions are taken to be of degree $p=3$ in $x$-direction, time-step size is chosen as $\Dt = 0.05$, and we use 1000 particles per cell. The tolerance of all appearing matrix solvers is set to $10^{-12}$.

\subsubsection{Cartesian Coordinates}

\begin{figure}[h!]
	\centering
	\includegraphics[width=\linewidth]{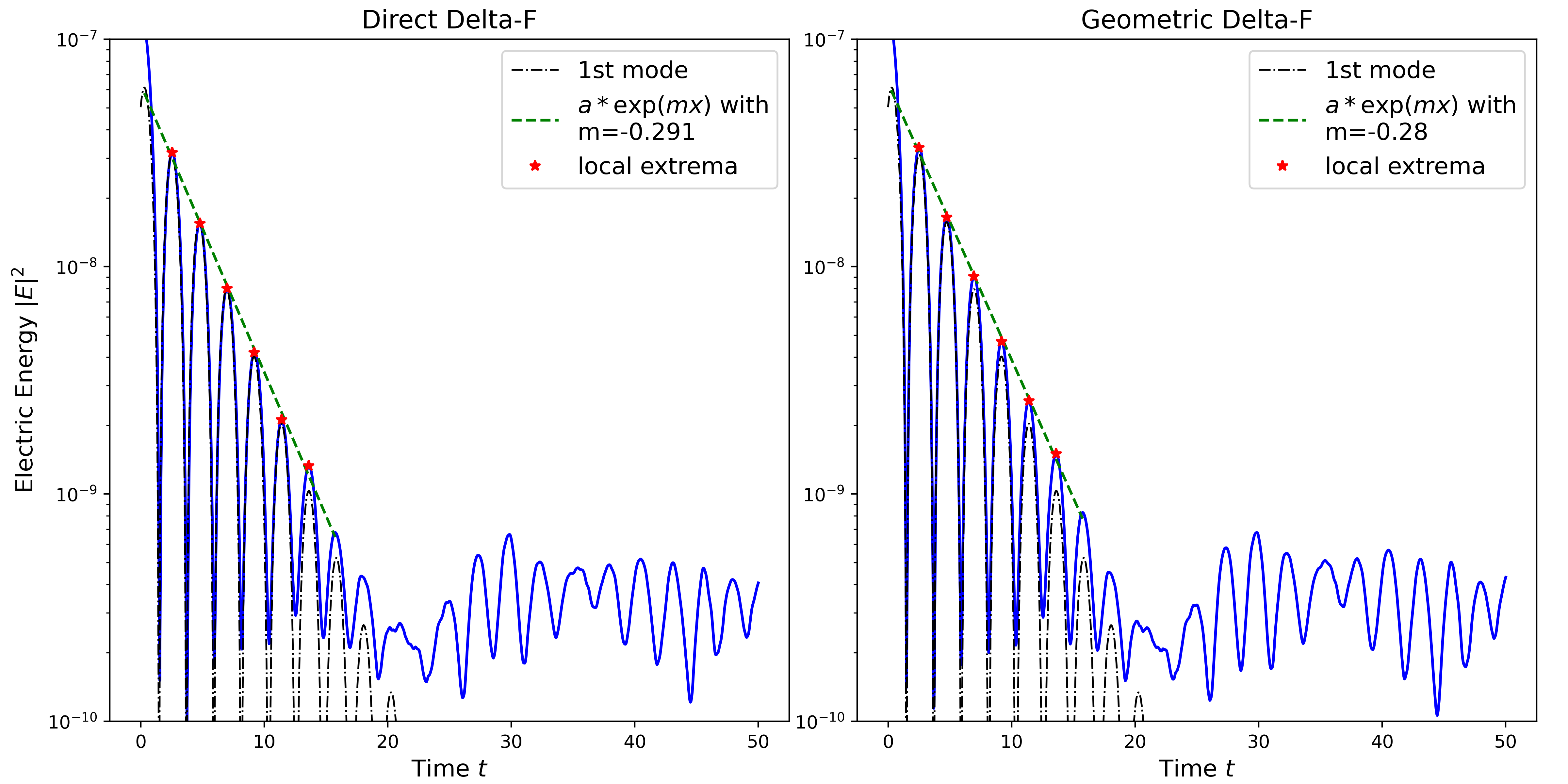}
	\caption{Weak Landau damping of the electric field energy, including growth rates for the damping phase. Our model (right) reproduces the analytical damping rate of $m = - 0.3066$ and reaches the same damping level as the comparison model (left) until saturation. For comparison, the first (dominant) mode of the exact solution (black dashed line) is also shown.}
	\label{fig:weaklandau}
\end{figure}

As a first check, we compare the models in the damping region (c.f. fig. \ref{fig:weaklandau}). The figure shows the electric field energy as a function of time until $t_{\text{end}} = 30$ for both models and of the dominant first mode of the exact solution \cite{Sonnendrücker_2015}. Additionally, we fit a slope to the first six local maxima. We see that both models give the same damping rate of about $m=-0.3 \pm 0.08$ and also saturate the damping between $10^{-10}$ and $10^{-9}$. In this comparison, the models produce consistent results.

\begin{figure}[h!]
	\centering
	\includegraphics[width=\linewidth]{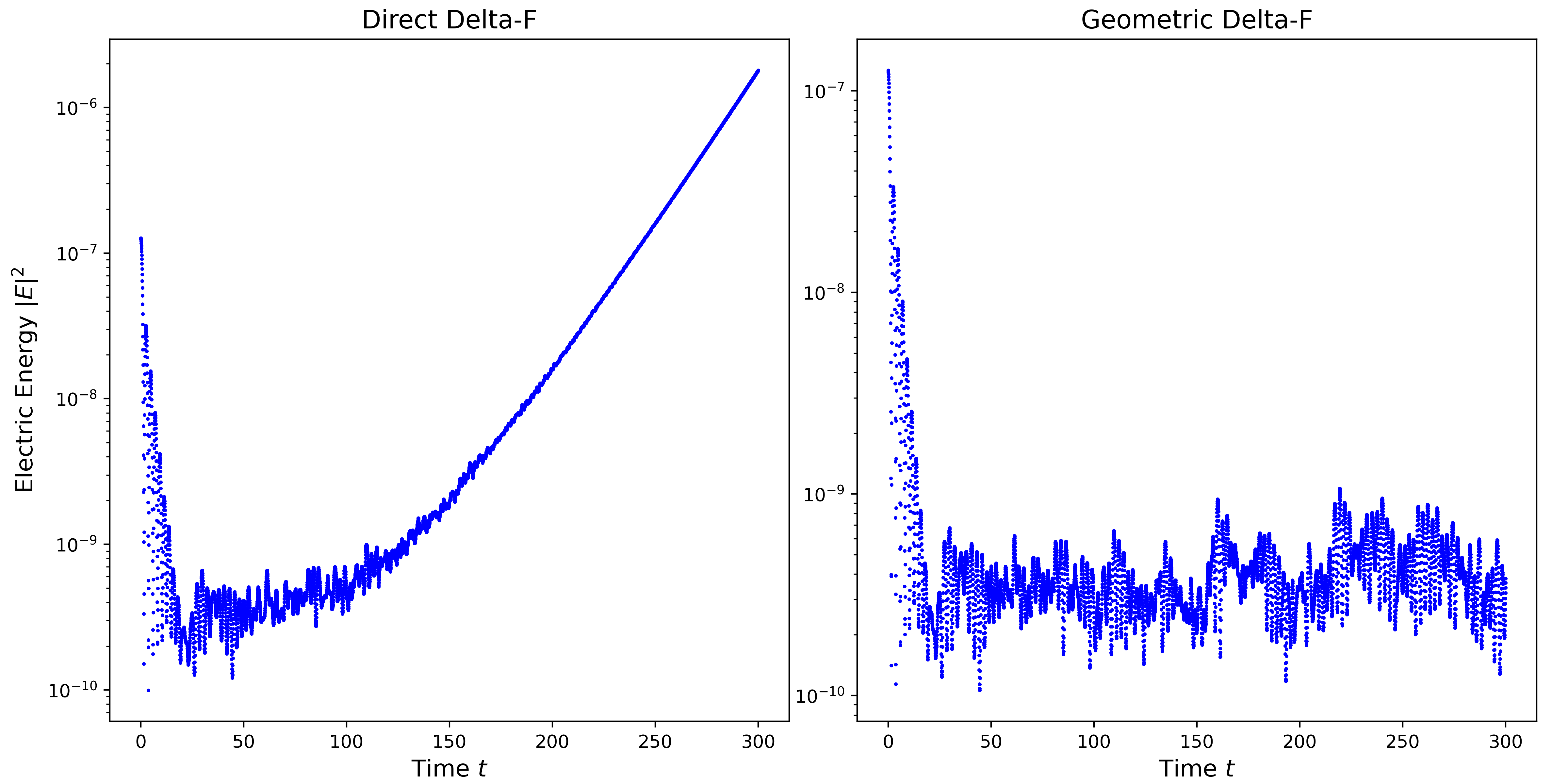}
	\caption{Long-time behaviour of the electric field energy in weak Landau damping. Our geometric model (right) keeps the amplitude of the energy at damped levels for long times, while the comparison model (left) shows an unphysical growth due to numerical noise after some time.}
	\label{fig:weaklandaulong}
\end{figure}

The advantage of the geometric formulation of our model becomes clear when we look at long-time simulations (c.f. fig. \ref{fig:weaklandaulong}). Again, the electric field energy is shown as a function of time, this time until $t_{\text{end}} = 200$. We see that our model keeps the level of the electric field energy at the same level of numerical accuracy. The comparison model shows a growth in electric field energy at late times, which is purely due to numerical noise and unphysical.

\begin{figure}[h!]
	\centering
	\includegraphics[width=\linewidth]{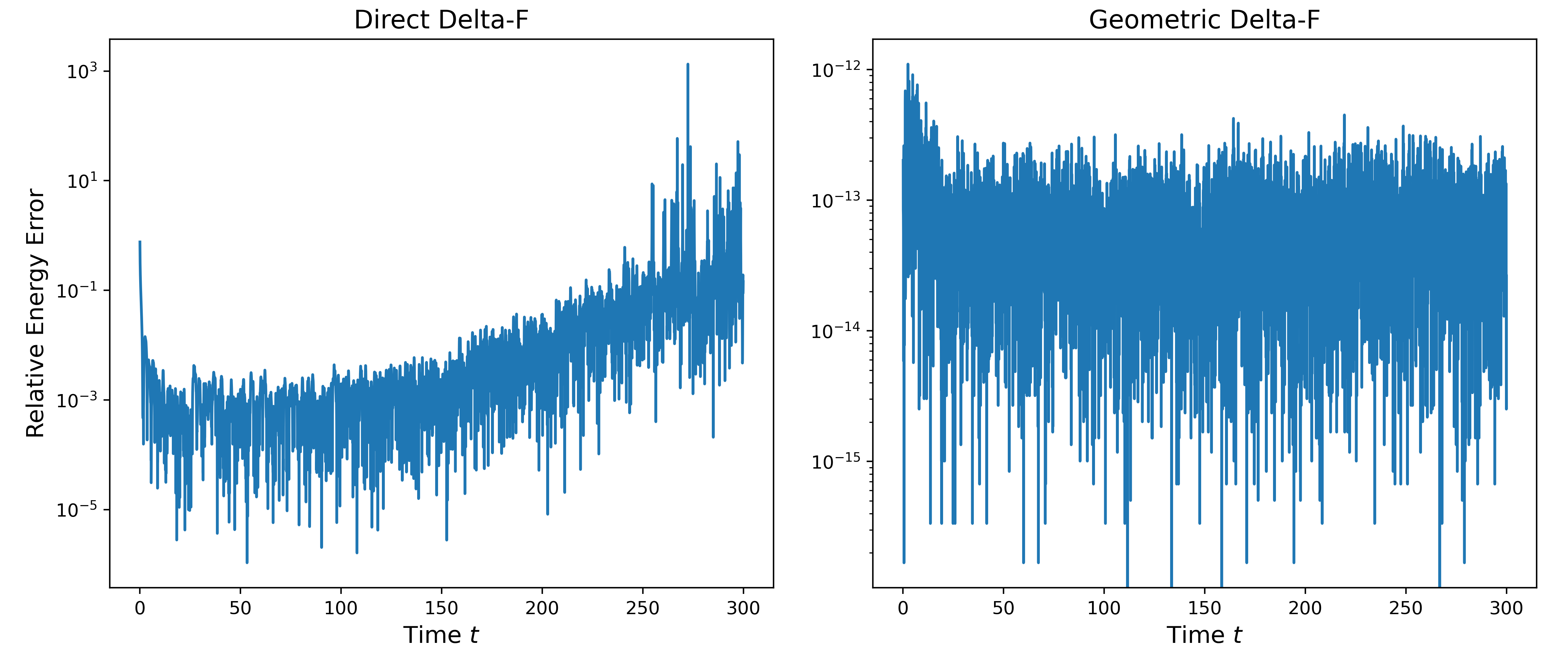}
	\caption{Relative error of the total energy \eqref{relative-energy-error} during weak Landau damping. Our model (right) keeps the error at the accuracy of the matrix inversion solver, while the comparison model (left) has a much larger error in total energy.}
	\label{fig:weaklandauentotrelerr}
\end{figure}

Furthermore, in fig. \ref{fig:weaklandauentotrelerr} we plot the absolute value of the relative error of the total energy
\begin{equation}\label{relative-energy-error}
	\abs{\frac{\abs{H_{\text{tot}}^{n+1} - H_{\text{tot}}^n}}{H_{\text{tot}}^n}} \, , \qquad \text{where $n$ is the time-step index}
\end{equation}
We take the absolute value of this whole expression since the control variate method also produces negative energies. We see that our geometric delta-$f$ formulation keeps the relative error at a level equivalent to the accuracy of the matrix solver used in the coupling-substep, sec. \eqref{sec:subsystem-two}.

\begin{figure}
	\centering
	\includegraphics[width=\linewidth]{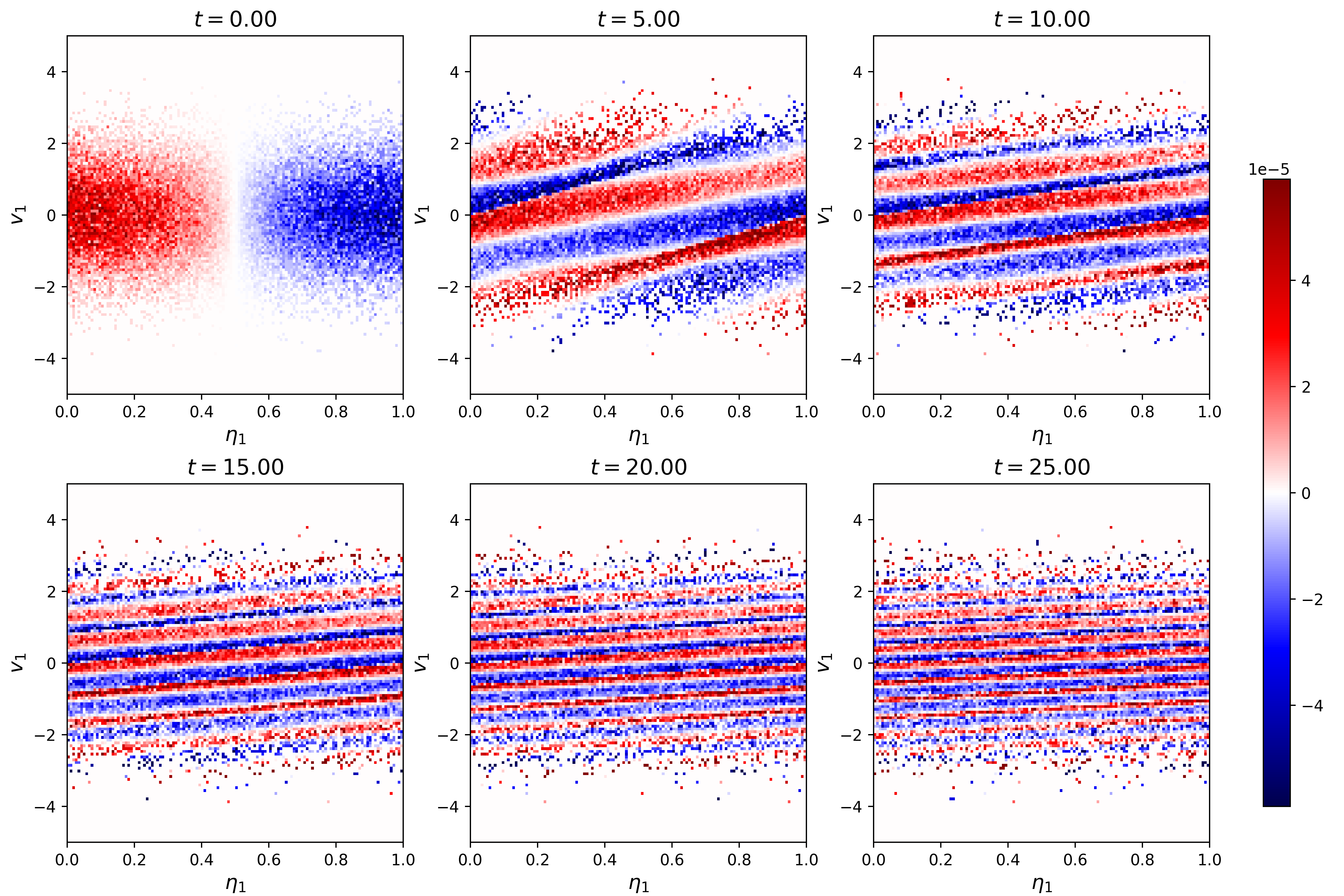}
	\caption{The evolution of the distribution function in phase space for different times during weak Landau damping. The filamentation can be clearly seen.}
	\label{fig:LVA_landau_phase_space}
\end{figure}

Finally, we plot the perturbation of the distribution function $f_1$ in phase space $v_1 - \eta_1$ for different times during the weak Landau damping (c.f. fig. \eqref{fig:LVA_landau_phase_space}). The initial condition is transported with different velocities, and the resulting filamentation can be seen very clearly. The cost of simulations for both models is similar, with our model taking 114$s$ to complete the long-time simulation, whereas the control-variate method takes 126$s$ on a M1 MacBook Pro run with 6 MPI processes.

\begin{figure}
	\centering
	\includegraphics[width=\linewidth]{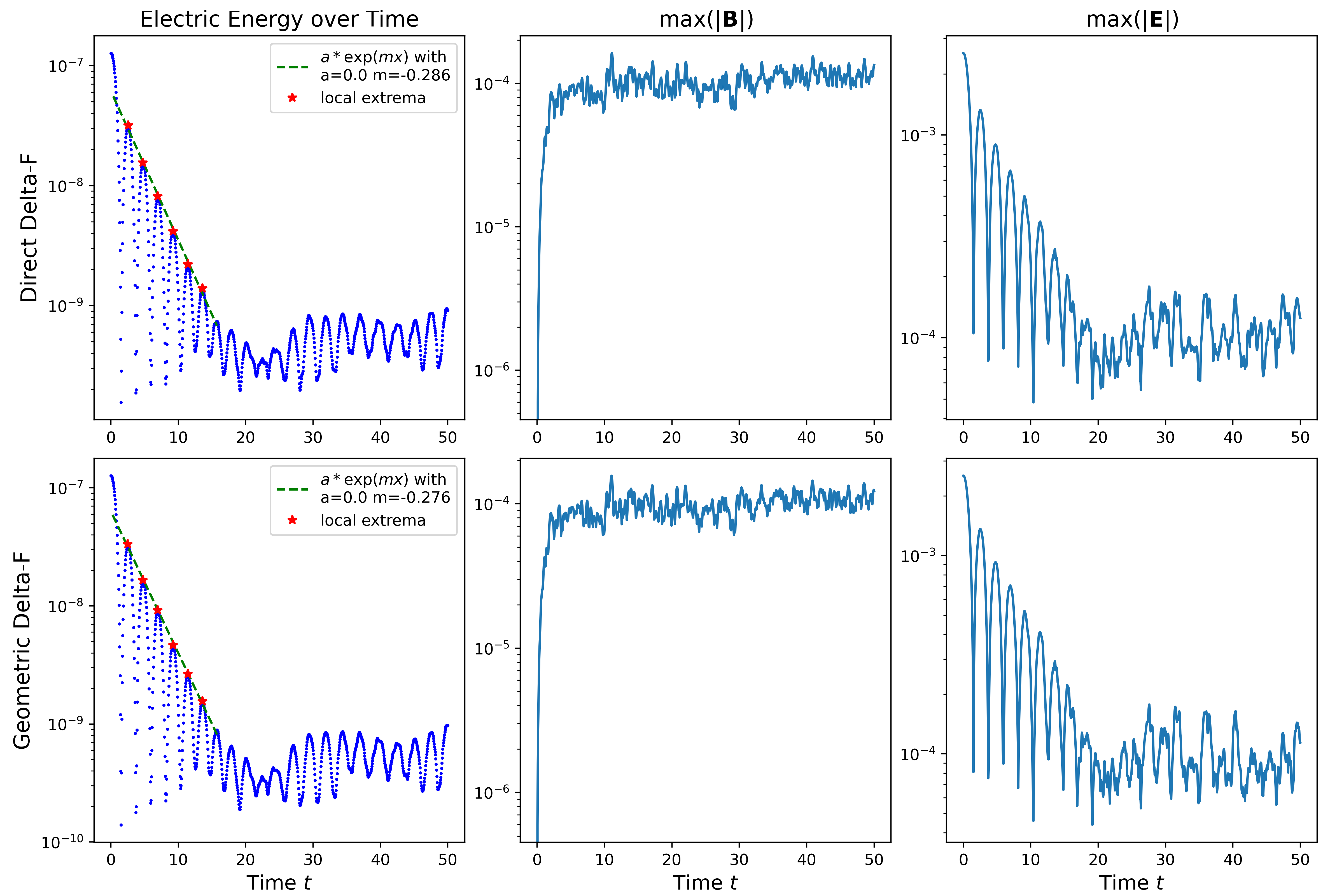}
	\caption{The weak Landau damping in Cartesian coordinates using our model with a dynamic magnetic field (bottom). We observe that our model gives a consistent damping rate of $m = -0.3066$ as the Vlasov--Ampère models in fig. \eqref{fig:weaklandau}. The magnetic field (center) is not identically zero, as is expected from theory, but rather has a small amplitude of $10^{-4}$ which is about the amplitude of the electric field (right) in the saturated regime, so it is due to numerical noise. This phenomenon also occurs in the comparison model (top).}
	\label{fig:landau_lvm}
\end{figure}

In fig. \eqref{fig:landau_lvm}, we also present the results using the linearized Vlasov--Maxwell model, in which the magnetic field is initialized as zero and should remain at zero over time. We first observe that the model reproduces the damping rate correctly. The magnetic field does not remain zero over time, but the amplitude of the magnetic field is of the same order as the electric field in the saturated phase, i.e., of the level of numerical noise. This level can be brought down by setting finer solver tolerances and increasing the number of particles for the simulation. Additionally, the divergence of the magnetic field is kept at identical zero throughout the whole simulation.

\subsubsection{Colella Coordinates}

\begin{figure}
	\centering
	\includegraphics[width=\linewidth]{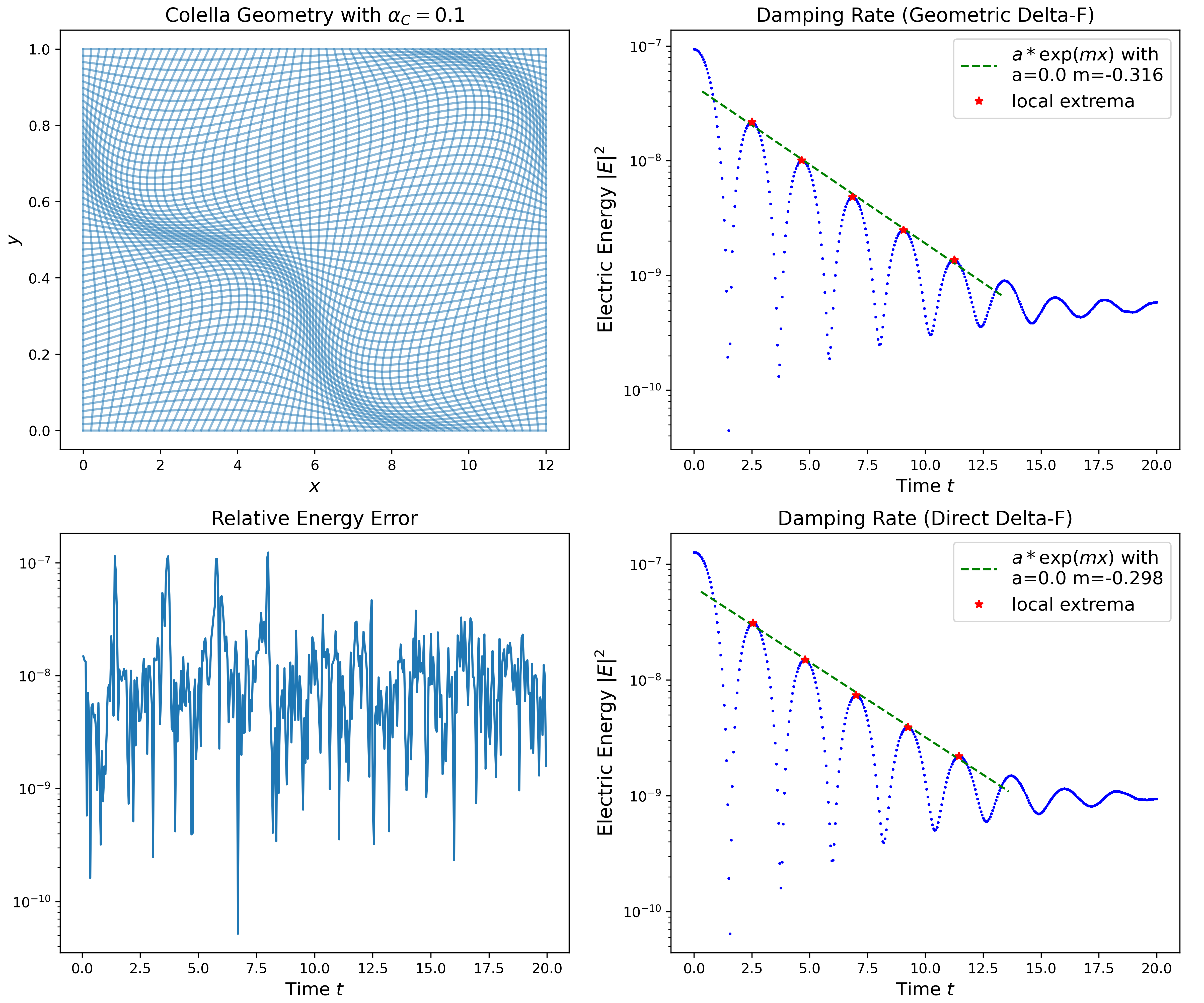}
	\caption{Weak Landau damping in curvilinear coordinates: Isolines of the Colella mapping, the damping rate of the linearized Vlasov--Ampère model (upper right) and the comparison model (lower right), and the relative energy error over time for our model (lower left), which is kept at the accuracy of the matrix solvers.}
	\label{fig:lva_landau_colella}
\end{figure}

Finally, for the weak Landau damping test case, we present results in the Colella mapping to verify that our model also works in curvilinear coordinates. The Colella mapping reads
\begin{equation}
	L_{\text{Colella}} : \hat{\Omega} \to \Omega \, , \quad \bh \mapsto \begin{pmatrix}
		L_x \left[ \eta_1 + \alpha_C \sin(2 \pi \eta_1) \sin(2\pi \eta_2) \right] \\
		L_y \left[ \eta_2 + \alpha_C \sin(2 \pi \eta_1) \sin(2\pi \eta_2) \right] \\
		L_z \eta_3
	\end{pmatrix} = \bx \, ,
\end{equation}
and mixes the first and second coordinates. We are essentially adding artificial curvature and the physical domain is still Cartesian. Boundary conditions are taken to be periodic in all directions and we choose $\alpha_C = 0.1$ for the appearing parameter. We consider it on a domain $\Omega = [0, 12] \times [0,1]^2$ with 24 cells in both $x$- and $y$-direction with degree 3 FE basis functions. Other parameters are as in the Cartesian case above. The initial condition is also the same with the perturbation going in $x$-direction. In fig. \eqref{fig:lva_landau_colella} we show the isolines of this geometry. We also show the damping rate of the model, which agrees with the comparison models. Furthermore, the relative energy error is kept small.

\subsection{Bernstein Waves}\label{sec:Bernstein-waves}

Bernstein waves \cite{Bernstein_waves} exist in plasmas subject to an external magnetic field. In contrast to the Landau waves discussed above, these modes are undamped and exist irrespective of the strength of the magnetic field. This discontinuous jump in existence of modes between vanishing and a small background magnetic field is known as the Bernstein--Landau paradox. Because the wavelengths of the Bernstein waves are of the order of the gyro-radius, they can only be captured by kinetic models, thus presenting a perfect candidate for another test case.

\subsubsection{Setup}

We take the domain to be $\Omega = [0,144] \times [-0.4, 0.4]^2$ and simulate until $T=2000$ with a time-step size of $\Delta t = 0.25$. The number of cells is chosen to be $1024 \times 1 \times 1$ with degree 3 splines in the first direction. The background magnetic field is taken as the unit vector in $z$-direction: $\bB_0 = e_z$. We take the background Maxwellian \eqref{f0-is-maxwellian} with a thermal velocity of $\vth = 0.2$ and a constant density $n_0(\bx) = 1$ such that no background electric field exists $\bE_0 = 0$. Since we want to excite all modes of the spectrum, we initialize the distribution function with random noise. The number for particles per cell is chosen to be 100. We discuss different initial conditions and their effect on the result in appendix \eqref{app:bernstein}.

\subsubsection{Results}

\begin{figure}[h!]
	\centering
	\includegraphics[width=0.95\linewidth]{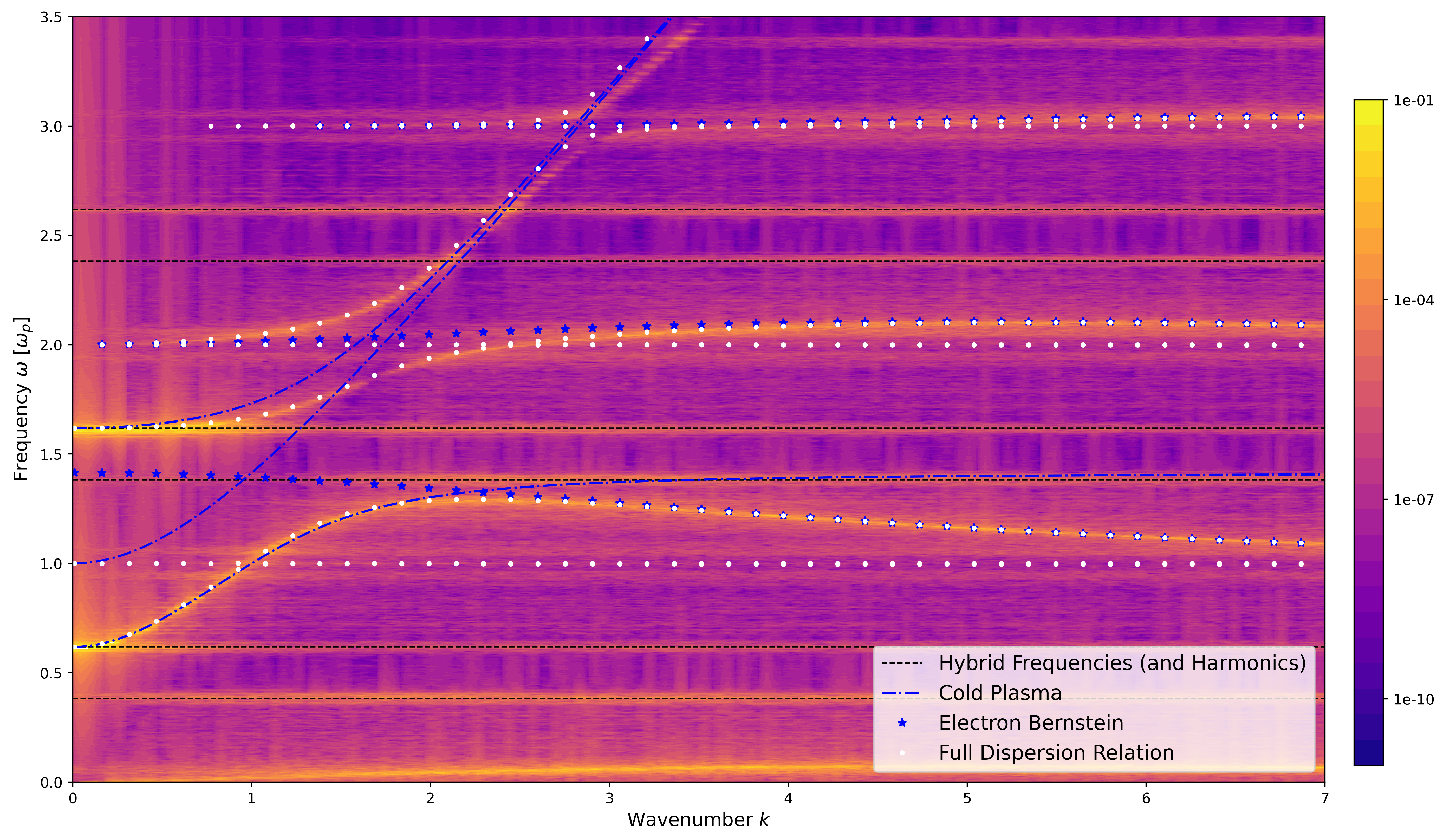}
	\caption{The perpendicular dispersion relation under a constant magnetic background field $\bB_0 = \be_z$. The Fourier transformed $x$-component of the electric field and the exact dispersion relation (blue, dotted), along with the electron Bernstein waves (blue, starred) and the cold plasma waves (blue, dot-dashed) are shown. The model reproduces the analytical dispersion relation very well, but also exhibits discrete Case--van Kampen modes at the harmonics of the hybrid frequencies (black, dashed).}
	\label{fig:lvmbernstein}
\end{figure}

To see the dispersion relation of the different waves, we Fourier transform the $x$-component of the electric field using the $N$-$d$ discrete Fourier transform function \texttt{scipy.fft.fftn}. Results are shown in fig. \eqref{fig:lvmbernstein} together with the exact dispersion relation \cite{Fitzpatrick_book}, the electron Bernstein waves, and the O- and X-mode of the cold plasma \cite{Meng_2025}, which are the asymptotics of the full dispersion relation. We see that the model reproduces the full dispersion relation very well in the whole spectrum. When comparing to the literature \cite{Yingzhe_2019, Xiao_2015_Bernstein, Kilian_2017, Xiao_2013_bernstein}, we also observe unphysical modes with frequencies not dependent on $k$ and have dominant modes located at harmonics of the hybrid frequencies (around half-integer multiples of $\omega_p$) and modes with lower amplitudes towards the Bernstein waves. They become very pronounced when using a sum of sines as the initial condition. These modes are discrete Case--van Kampen waves whose spectrum is continuous for the continuous model but becomes discretized by the discretization. We study them and discuss properties in app. \eqref{app:bernstein}. There is also a wave closely above $\omega = 0$ which vanishes as the time-step size $\Delta t \to 0$.

% !TeX spellcheck = en_GB

\section{Summary \& Outlook}\label{sec:summary_outlook}

In this work, we presented a Hamiltonian framework of the linearized Vlasov--Maxwell model with a Maxwellian background distribution function $f_0$ and consistent electromagnetic background fields $\bE_0$ and $\bB_0$. The system consists of the Maxwell equations \eqref{final-cart-maxwell-equations} and a Vlasov-type equation \eqref{final-cart-vlasov-equation}, which has a source term proportional to the electric field perturbation. The Hamiltonian of the model \eqref{final-cart-Hamiltonian} is quadratic in all three variables, which is also to be expected when deriving it from the energy--Casimir method \cite{Holm1985, Morrison_Ideal_Fluid, Shadwick_PhDThesis}. Furthermore, a non-canonical Poisson bracket was found that generates the equations of motion, completing the Hamiltonian system. The equations were then formulated in curvilinear coordinates using a diffeomorphism $L$ as a mapping \eqref{mapping} to enable using the model for simulations of realistic tokamak or stellerator geometries.

The discretization was done using the GEMPIC framework, in which the electromagnetic fields are understood as coefficient functions of forms living in an exact de Rham sequence. This FEEC approach lets us keep the geometric properties of the theory at the discrete level. The distribution function was discretized using a particle-in-cell method, for which we developed the tools necessary to handle the source term. We also discretized the Hamiltonian and the Poisson bracket. To integrate the system in time, we employed a Poisson splitting and found numerical methods to solve each of the resulting subsystems. The implementation was done in \texttt{STRUPHY} \cite{StruphyPaper, StruphyDoc}.

The model was tested on two test cases: firstly, the Landau damping, both in 1D Cartesian coordinates and in a 2D curved geometry, where we compared it to the Vlasov--Maxwell system with a control variate method and the linearized Vlasov--Maxwell model using the direct delta-f method. Our model gave consistent results in the domain of validity with the comparison models, i.e., a correct damping rate and final magnitude of the electric field energy. But in long-time simulations, it is where our geometric methods trump the simple noise reduction technique (fig. \ref{fig:weaklandaulong} and fig. \ref{fig:lva_landau_colella}): the amplitude of the electric field energy is kept at a low value while the numerical errors accumulate in the comparison models and blow up eventually. Also, energy conservation (fig. \ref{fig:weaklandauentotrelerr}) is constantly respected (up to solver accuracy). The second test case was the dispersion relation under a constant magnetic background field (fig. \eqref{fig:lvmbernstein}), where our model correctly reproduced the analytical dispersion relation, in particular the electron Bernstein waves for larger wave-numbers $k$. We observed other, unphysical modes at multiples of the hybrid frequencies which were found to be discrete Case--van Kampen stemming from the discretization.

One obvious follow-up work is the generalization of this model to more background distribution functions. The model now does not allow for shifts or anisotropies in the velocity part (c.f. eq. \eqref{f0-is-maxwellian}). A first step would be the introduction of shifts and different thermal velocities in all directions. However, the constructive way of the energy--Casimir method does not work with more general backgrounds, for example, a shift in the background distribution function and vanishing background fields
\begin{align}
	f_0 & = \frac{1}{\sqrt{(2 \pi)^3}} \exp\left( - \frac{\abs{\bv - \bu}^2}{2} \right) \, , & \bE_0 & = 0 \, , & \bB_0 & = 0 \, ,
\end{align}
would need a linear term in $\bv$ in the functional, namely $\cF \sim \bu \cdot \bv \, f$, but this is incompatible with $\bB_0 = 0$. The method of dynamically accessible variables \cite{Morrison_1989, Morrison_1990} could overcome this problem.

\textbf{Acknowledgements} The authors would like to thank Li Yingzhe, Klaus Hallatschek, and Philip J. Morrison for insightful discussions.

\appendix
% !TeX spellcheck = en_GB

\section{Analysis of Bernstein Waves}\label{app:bernstein}

Here we supply some additional material discussing the initial conditions, discrete Case-van Kampen modes, and providing more results obtained in the simulating of Bernstein waves (c.f. sec. \eqref{sec:Bernstein-waves}). For better illustration purposes, this section uses a slightly colder plasma with $\vth=0.06$ than in the main text.

\subsection{Different Initial Conditions}

\begin{figure}[h!]
	\centering
	\includegraphics[width=0.9\linewidth]{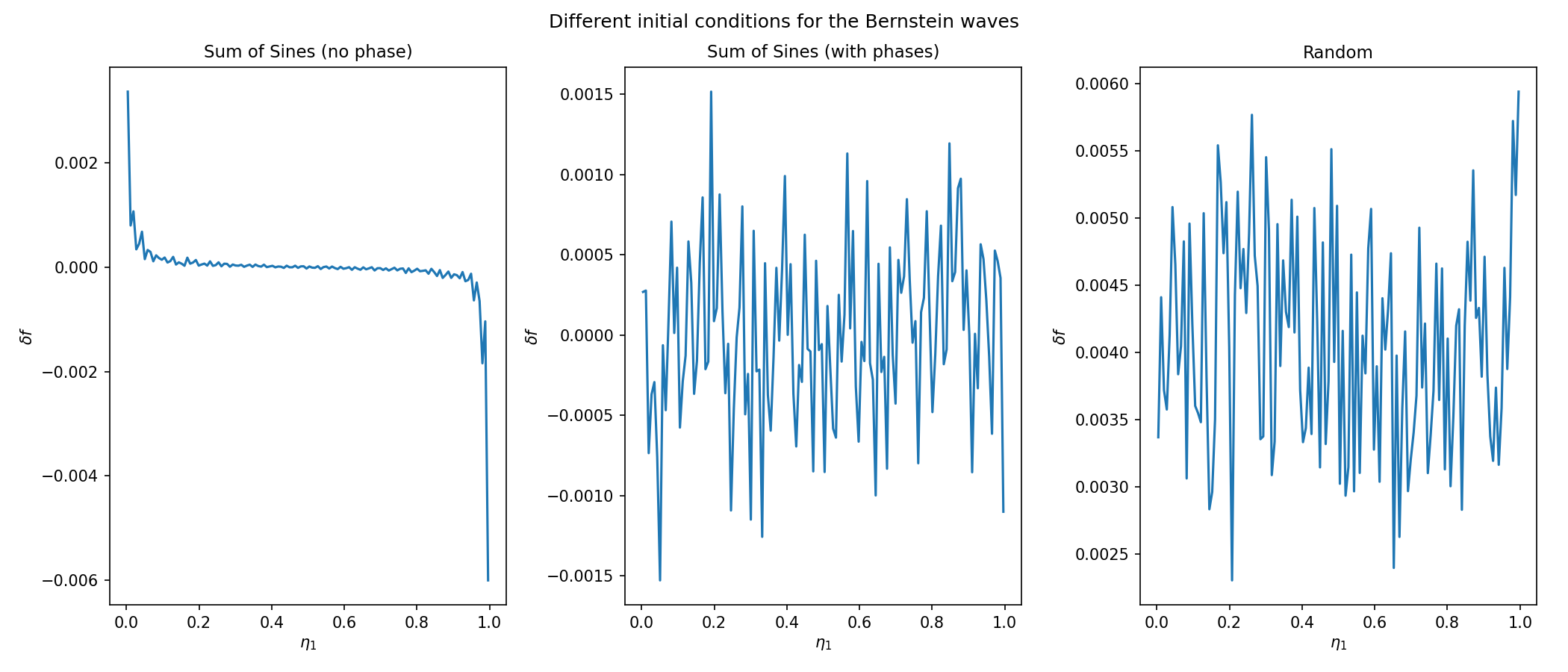}
	\caption{Different initial conditions for the simulation of Bernstein waves \eqref{sec:Bernstein-waves}. Shown are a sum of sines \eqref{bernstein-sum-of-sines-init} without phases (left) and with phases proportional to $k^2$ (center). We see that the inclusion of these phases smoothens the initial condition out quite a bit, such that it resembles the random noise (right).}
	\label{fig:lvmbernsteindifferentinits}
\end{figure}

\begin{figure}
	\centering
	\includegraphics[width=0.95\linewidth]{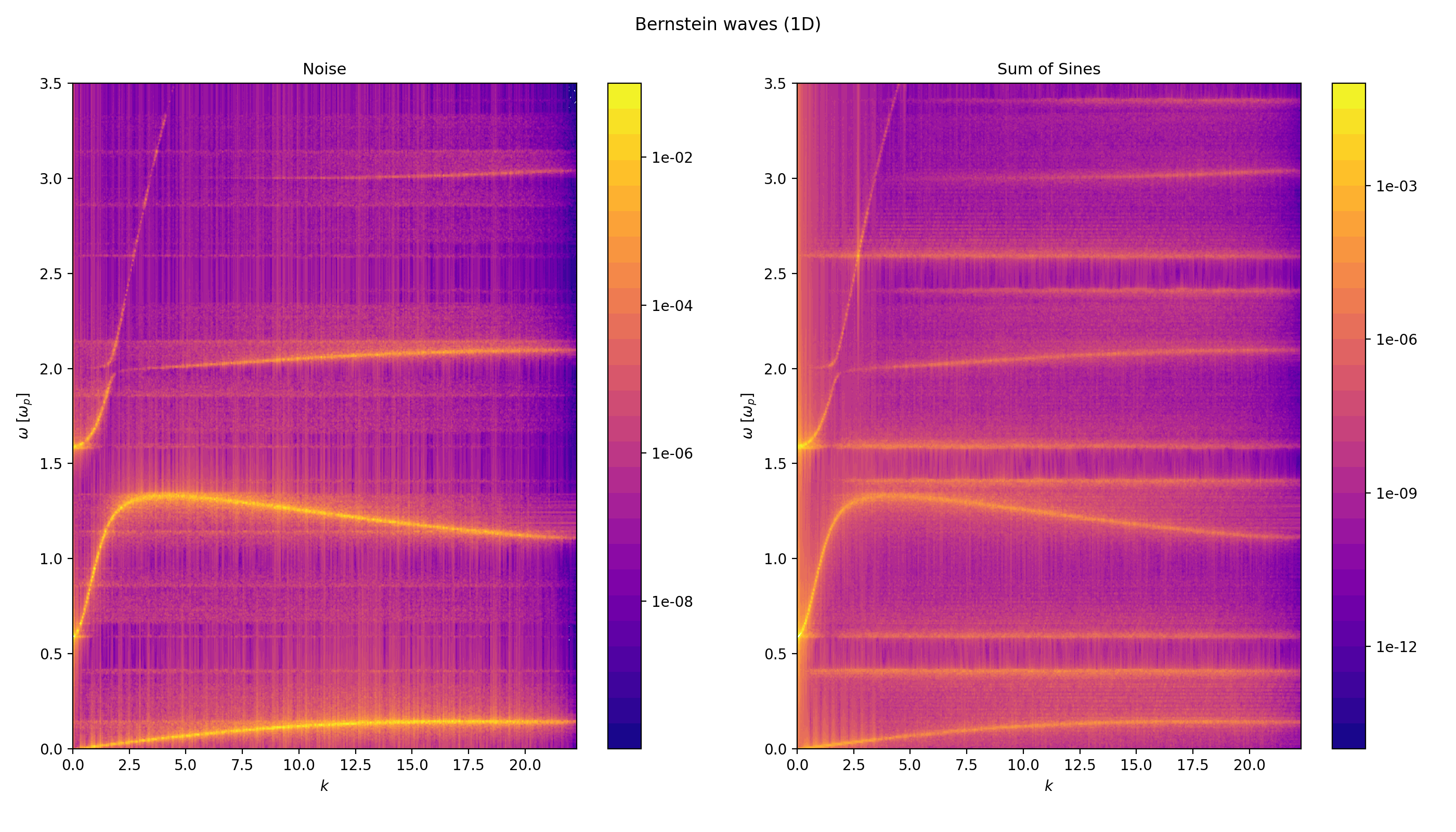}
	\caption{A comparison of the Fourier transformed $x$-component of the electric field for the two initial conditions of random noise and sum of sines with phases \eqref{bernstein-sum-of-sines-init}. We observe that the horizontal lines around half-integer $\omega_p$ become very pronounced with an amplitude comparable to the physical waves.}
	\label{fig:lvmbernsteincomparison}
\end{figure}

We first show the difference between the sum of sines
\begin{equation}\label{bernstein-sum-of-sines-init}
	f_1(t=0, \bx, \bv) = \alpha \, f_0(\bv) \, \sum_{k=0}^{80} \sin \left( 2 \pi \frac{k x}{L_x} + 2 \pi \sqrt{2} k^2 \right) \, , \qquad \text{where } \alpha = 10^{-4} \, ,
\end{equation}
with and without phases, and a (pseudo-) random drawing of marker weights in figure \eqref{fig:lvmbernsteindifferentinits}. Although the added phases do remove the peaked character of the pure sum of sines, and looks similar to the purely random noise initial condition, it exhibits the same modes as the peaked distribution \eqref{fig:lvmbernsteincomparison}.

\subsection{Discrete Case--van Kampen Modes}

\begin{figure}
	\centering
	\includegraphics[width=0.95\linewidth]{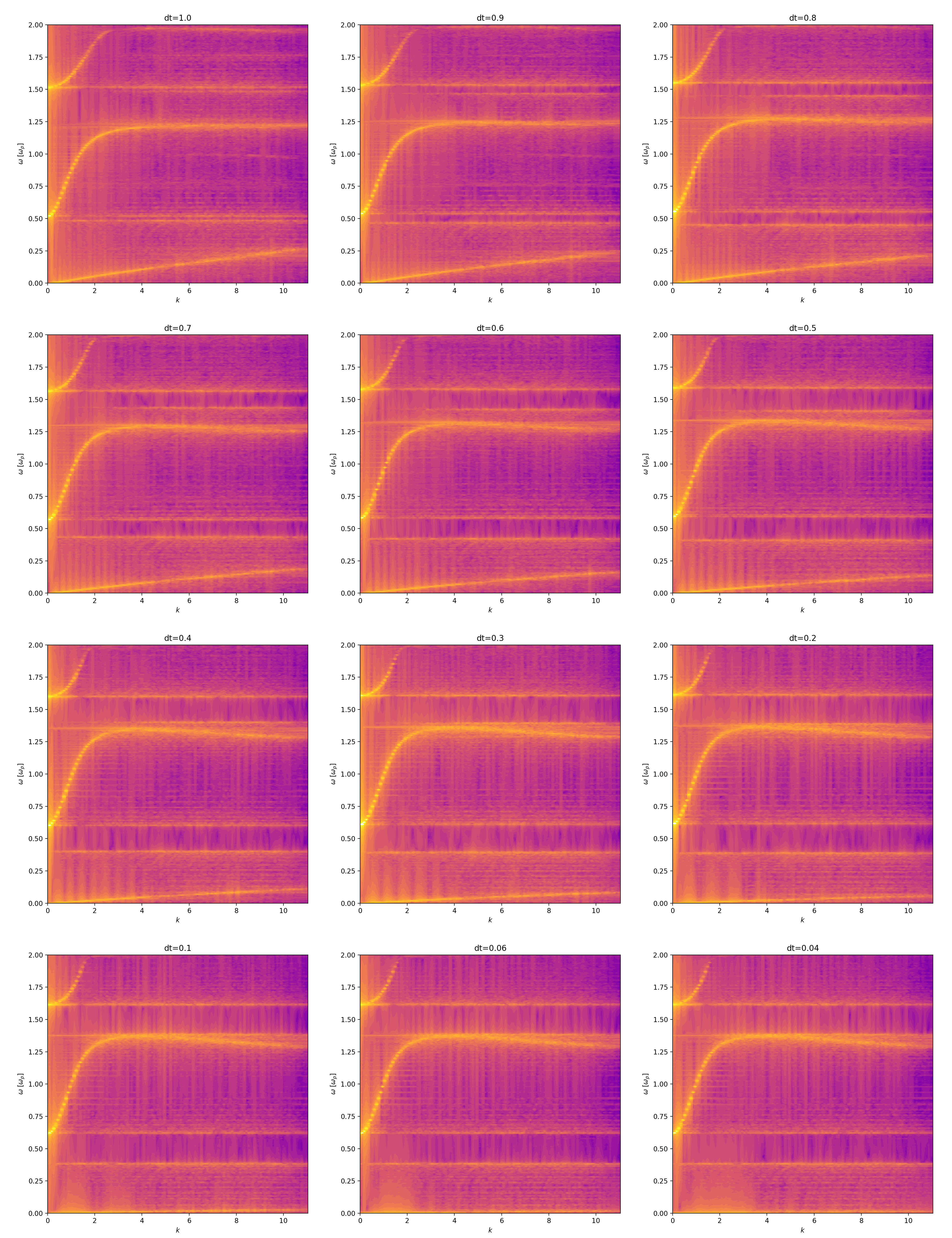}
	\caption{The dispersion plot for different time-step sizes. Shown is the Fourier transformed $x$-component of the electric field for the sum of sines initial condition. It can be seen that the "bands" of the Case--van Kampen modes open for smaller $\Delta t$.}
	\label{fig:lvmbernsteindtdispersions}
\end{figure}

\begin{figure}
	\centering
	\includegraphics[width=0.95\linewidth]{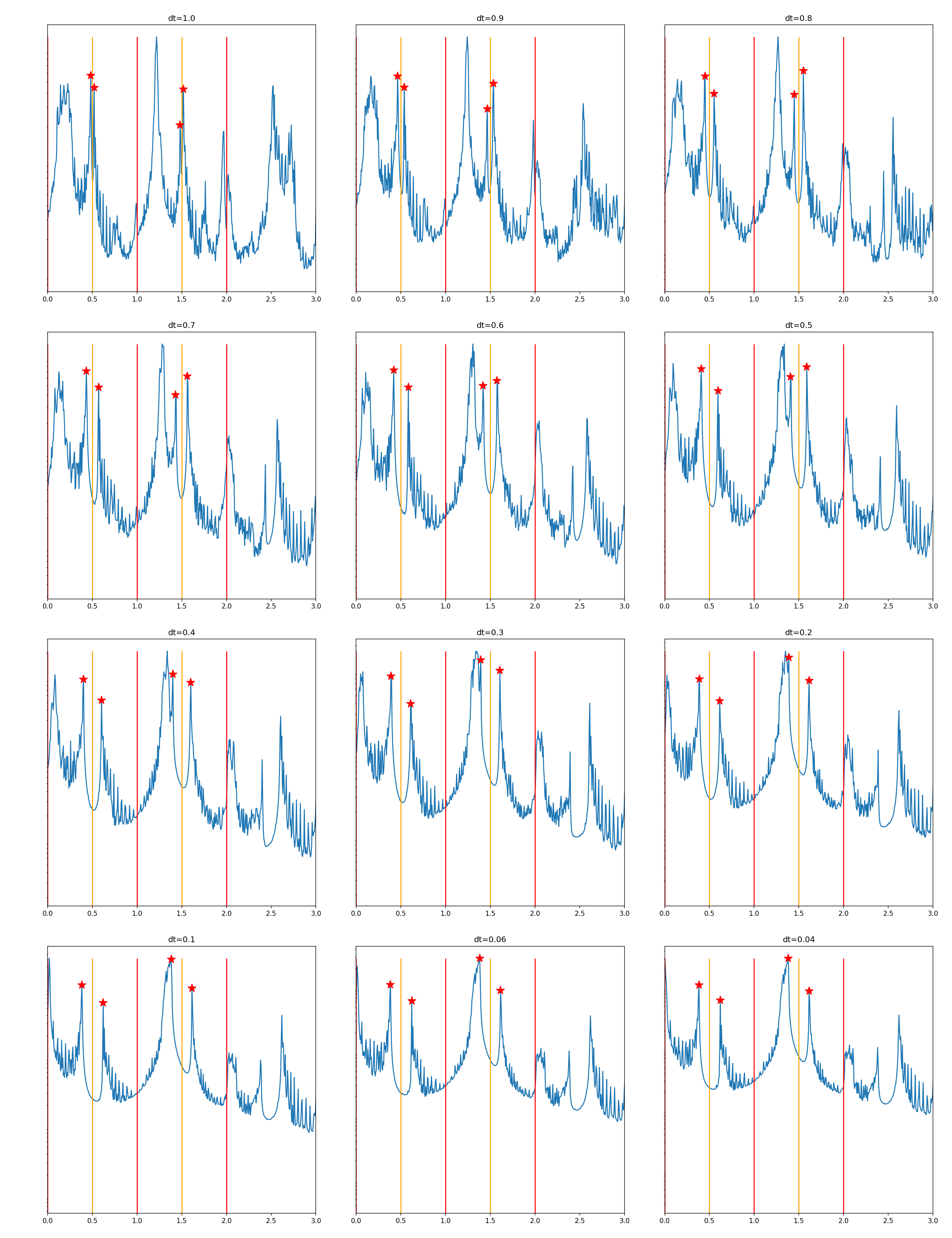}
	\caption{Shown is the averaged Fourier transformed $x$-component of the electric field for the sum of sines initial condition. The average is taken from $k_\text{start}=4$ until the end of the $k$-domain. The peaks (highlighted with red stars) around the half-integer multiples of the plasma frequency (vertical orange lines) are moving apart as $\Delta t$ decreases.}
	\label{fig:lvmbernsteindtdistances}
\end{figure}

The horizontal lines that can be seen in the Fourier transformed electric field are discrete Case--van Kampen modes \cite{van_Kampen_1955, Case_1959, Case_1967_book}. These waves have a continuous spectrum but are discretized due to the discretization of the equations of motion.  We show the dispersion relations for different time-step sizes (fig. \eqref{fig:lvmbernsteindtdispersions}) and the averaged amplitude starting at $k=4$ (fig. \eqref{fig:lvmbernsteindtdistances}). These simulations were run with $[256, 1, 1]$ cells on a domain $[0, 72] \times [-0.4, 0.4]^2$, and with 1000 particles per cell.

\begin{figure}[h!]
	\centering
	\includegraphics[width=0.55\linewidth]{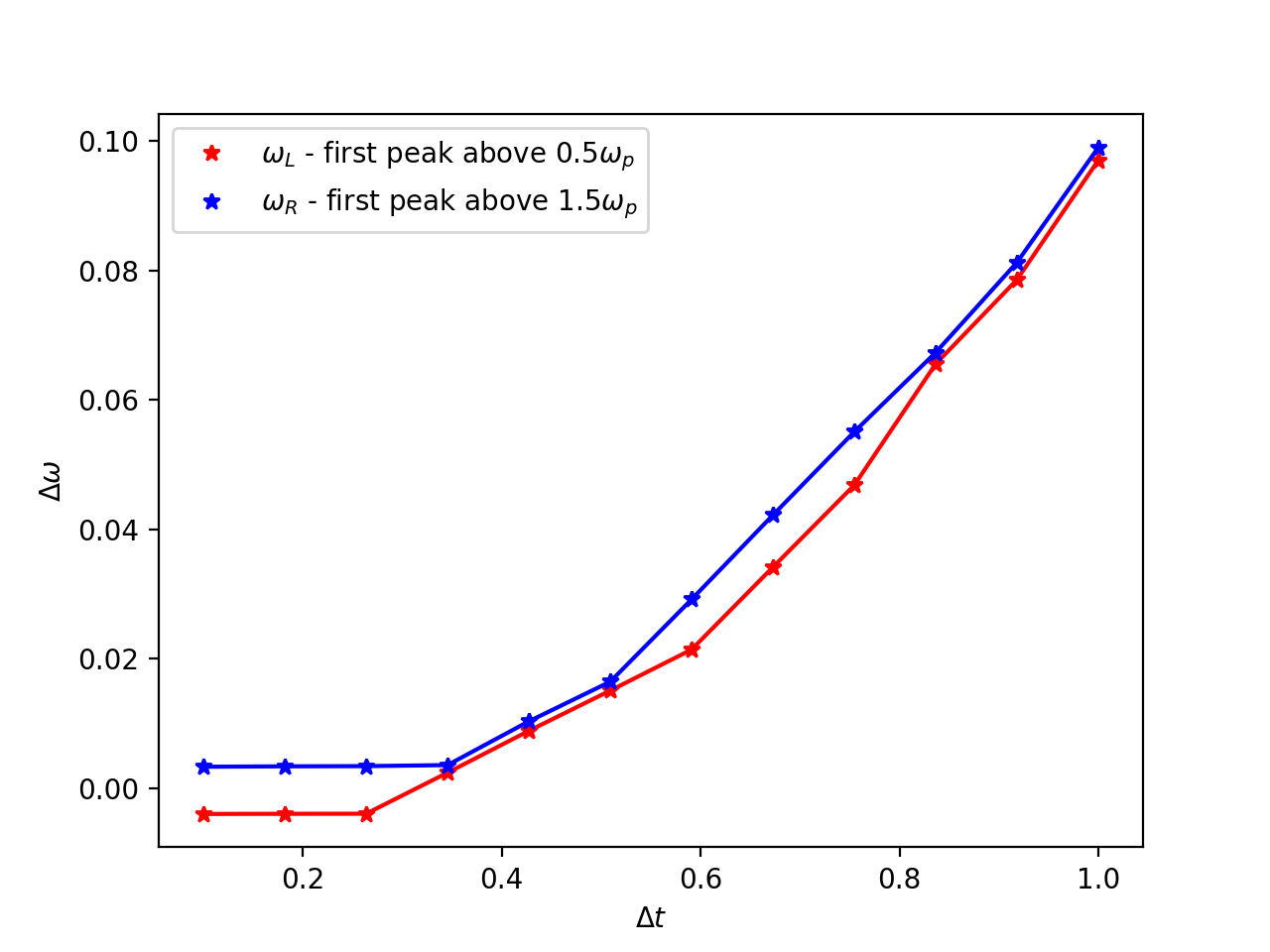}
	\caption{The difference between the analytical value of the abscissas of the slow and fast $X$-wave and the numerical value for different time-step sizes $\Delta t$.}
	\label{fig:lvmbernsteindtpeakplot}
\end{figure}

We see that the discrete Case--van Kampen appear symmetrically around half integer values of $\omega_p = 1 = \omega_c$. These "bands" are close together initially but move apart as $\Delta t$ decreases and the temporal discretization error becomes smaller. The upper edge of lowest band is always marks the abscissas of the slow $X$-wave which moves closer to its analytical values of $\omega_L$; the same happens to the fast $X$-wave and its abscissas $\omega_R$:
\begin{subequations}
	\begin{align}
		\omega_L & = \frac{\omega_c}{2} \left( \sqrt{1 + 4 \left(\frac{\omega_p}{\omega_c}\right)^2} - 1 \right) = \frac{\sqrt{5} - 1}{2} \, , \\
		\omega_R & = \frac{\omega_c}{2} \left( \sqrt{1 + 4 \left(\frac{\omega_p}{\omega_c}\right)^2} + 1 \right) = \frac{\sqrt{5} + 1}{2} \, .
	\end{align}
\end{subequations}
Plotting the difference between the analytical and the numerical values for $\omega_L$ and $\omega_R$ respectively (c.f. fig \eqref{fig:lvmbernsteindtpeakplot}) we see that the difference approaches zero as $\Delta t \to 0$.

We also observe in fig. \eqref{fig:lvmbernsteindtdispersions} that the lowest (unphysical) mode which starts at $\omega = 0$ becomes flatter and eventually goes to zero as $\Delta t \to 0$. Particle number (per cell), number of cells, degree of splines, and strength of the background magnetic field were not seen to influence the number or behaviour of the Case--van Kampen modes.

\subsection{Additional Material}

\begin{figure}[h!]
	\centering
	\includegraphics[width=0.9\linewidth]{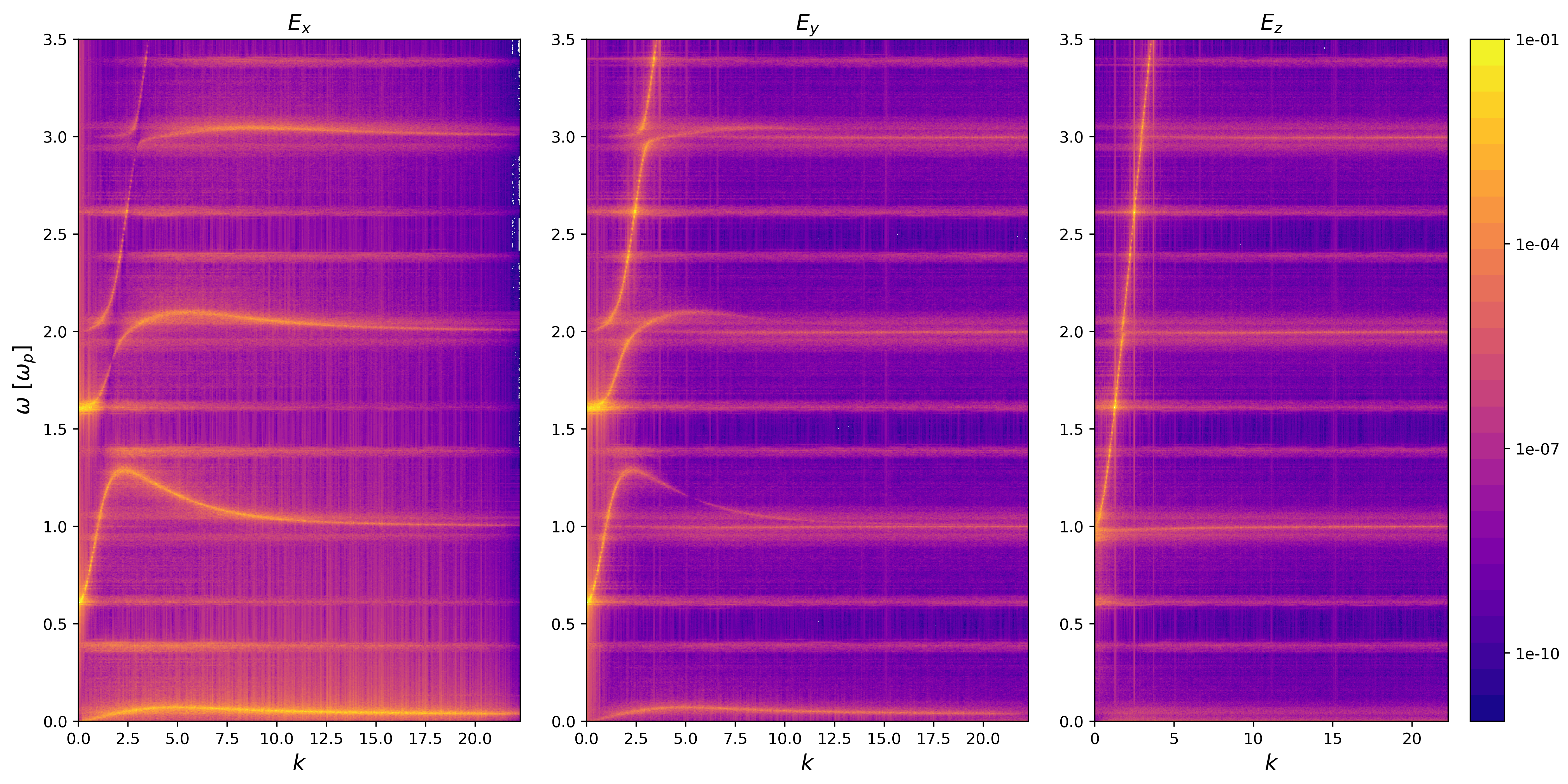}
	\caption{A comparison of the Fourier transformation of each component of the electric field. The Bernstein waves are only visible in the perpendicular directions while the electromagnetic and discrete Case--van Kampen modes are present in all components.}
	\label{fig:lvmbernsteincomponents}
\end{figure}

Additionally to the $x$-component of the electric field, which was discussed in the text, we show all the Fourier transformed components of the electric field in fig. \eqref{fig:lvmbernsteincomponents}. As expected, the Bernstein modes are only present in the components perpendicular to the external magnetic field and are absent in the direction parallel to it ($z$-component). There, the electromagnetic waves can also propagate freely without interacting with the Bernstein waves as in the perpendicular directions. The discrete Case--van Kampen modes are present in all components.

\printbibliography

\end{document}